\documentclass[12pt]{article}
\usepackage{amsmath}
\usepackage{latexsym}
\usepackage{amssymb}
\newtheorem{thm}{Theorem}[section]
\newtheorem{la}[thm]{Lemma}
\newtheorem{Defn}[thm]{Definition}
\newtheorem{Remark}[thm]{Remark}
\newtheorem{Note}[thm]{Note}

\newtheorem{Example}[thm]{Example}
\newtheorem{Examples}[thm]{Examples}
\newtheorem{Problems}[thm]{Problems}

\newtheorem{Problem}[thm]{Problem}
\newtheorem{Convention}[thm]{Convention}
\newtheorem{Number}[thm]{\!\!}
\newenvironment{defn}{\begin{Defn}\rm}{\end{Defn}}

\newenvironment{example}{\begin{Example}\rm}{\end{Example}}
\newenvironment{examples}{\begin{Examples}\rm}{\end{Examples}}

\newenvironment{rem}{\begin{Remark}\rm}{\end{Remark}}
\newenvironment{numba}{\begin{Number}\rm}{\end{Number}}
\newenvironment{proof}{{\noindent\bf Proof.}}%
                  {\nopagebreak\hspace*{\fill}$\Box$\medskip\medskip\par}   
\newcommand{\Punkt}{\nopagebreak\hspace*{\fill}$\Box$}
\newcommand{\wb}{\overline}
\newcommand{\ve}{\varepsilon}
\newcommand{\at}{\symbol{'100}}

\newcommand{\wt}{\widetilde}

\newcommand{\impl}{\Rightarrow}

\newcommand{\mto}{\mapsto}

\newcommand{\N}{{\mathbb N}}
\newcommand{\R}{{\mathbb R}}
\newcommand{\bO}{{\mathbb O}}
\newcommand{\F}{{\mathbb F}}
\newcommand{\K}{{\mathbb K}}

\newcommand{\Q}{{\mathbb Q}}

\newcommand{\Z}{{\mathbb Z}}

\newcommand{\C}{{\mathbb C}}

\newcommand{\wh}{\widehat}

\newcommand{\sub}{\subseteq}

\DeclareMathOperator{\im}{im}

\DeclareMathOperator{\pr}{pr}

\DeclareMathOperator{\id}{id}

\newcommand{\sbull}{{\scriptscriptstyle \bullet}}

\begin{document}
\renewcommand{\thefootnote}{\fnsymbol{footnote}}
\begin{center}
{\Large\bf Comparison of some notions of {\boldmath
$C^k$}-maps in\\[2.7mm]
multi-variable non-archimedian analysis}\\[6mm]
{\bf Helge Gl\"{o}ckner\footnote{\,These studies are part
of the project 436 RUS 17/67/05 by
the German Research Foundation (DFG).}}\vspace{2.5mm}
\end{center}
\renewcommand{\thefootnote}{\arabic{footnote}}
\setcounter{footnote}{0}
\begin{abstract}\vspace{1mm}
\hspace*{-7.2 mm}
Various definitions of $C^k$-maps
on open subsets
of finite-dimensional
vector spaces over a complete valued field
have been proposed in the literature.
We show that the $C^k$-maps
considered by Schikhof and De Smedt
coincide with those of Bertram, Gl\"{o}ckner
and Neeb.
By contrast, Ludkovsky's $C^k$-maps need not be $C^k$
in the former sense,
at least in positive characteristic.
We also compare
various types~of\linebreak
H\"{o}lder differentiable maps
on finite-dimensional
and metrizable spaces.
\vspace{2.6mm}
\end{abstract}
{\footnotesize {\em Classification}:
26E30 (primary); 26E20, 46A16, 46G05, 46S10\\[1mm] 
{\em Key words}: Non-archimedian analysis,
differentiable map,
partial difference quotient,\linebreak
directional difference quotient,
differentiability property, weak differentiability}\vspace{5mm}
\begin{center}
{\Large\bf Introduction}\vspace{1.6mm}
\end{center}
Various concepts of $C^k$-maps on subsets of
finite-dimensional vector spaces have been used
in the literature on non-archimedian analysis.
Schikhof's textbook
\cite{Sch} gave a comprehensive discussion
of the single-variable calculus
of $C^k$-maps over a complete ultrametric
field~$\K$, and suggested a
definition of multi-variable $C^k$-maps (in \S84),
which was then elaborated
by De~Smedt~\cite{DSm}.
Ludkovsky introduced a notion
of $C^k$-map between open
subsets of locally convex spaces
over a finite extension~$\K$ of $\Q_p$
(see \cite[Definition~2.3]{Ld1}
and \cite[Part~I, Definition~2.3]{Ld2}
for the case of Banach spaces,
\cite[Part~II, Remark~4.4]{Ld2} for the general case).
Recently, Bertram, Gl\"{o}ckner and Neeb~\cite{BGN}
introduced a notion of
$C^k$-map between open subsets
of arbitrary (Hausdorff) topological
vector spaces over a (non-discrete)
topological field~$\K$.
While the definition of $C^k$-maps by
Schikhof and De Smedt is based
on the existence of continuous
extensions to certain partial difference quotients,
the definition of Bertram et al.\
and Ludkovsky's definition
are based on continuous extendibility
of certain iterated directional difference
quotients.
The primary goal of this paper
is to compare these notions
of $C^k$-maps, and some related concepts.
To describe our main results,
let $E$ and $F$ be topological vector spaces
over a topological field~$\K$,
$k\in \N_0$,
and $f\colon U\to F$ be a map
on an open set $U\sub E$.
We start with a special case of Theorem~\ref{mainthm},
which generalizes a result for functions
of a single variable obtained in~\cite[Proposition~6.9]{BGN}.\\[2.5mm]
{\bf Theorem~A.}
\emph{If $E=\K^d$ for some $d\in \N$,
then $f$ is $C^k$ in the sense
of Bertram et al.\
if and only if $f$ is $C^k$ in the
sense of Schikhof and De Smedt.}\\[2.5mm]
If $\K$ is a valued field,
then variants of the two approaches just discussed
can be used to define $k$ times
H\"{o}lder differentiable maps
with H\"{o}lder exponent $\sigma\in \; ]0,1]$
($C^{k,\sigma}$-maps, for short).
As a special case of Theorem~\ref{mainthm},
we have:\\[2.5mm]
{\bf Theorem~B.}
\emph{If $E=\K^d$ for some $d\in \N$,
then $f$ is $C^{k,\sigma}$ in the sense
of Bertram et al.\
if and only if $f$ is $C^{k,\sigma}$ in the
sense of Schikhof and De Smedt.}\\[2.5mm]
By contrast, the mappings introduced
by S.\,V. Ludkovsky differ from
the preceding ones, if his definition
is used for fields of positive characteristic.
We show by example (see Theorem~\ref{countludk}):\\[2.5mm]
{\bf Theorem~C.}
\emph{For each local field $\K$ of positive
characteristic, there is a map $f\colon \bO\to \K$
on $\bO:=\{z\in \K\colon |z|\leq 1\}$
which is $C^\infty$ in Ludkovsky's
sense, but not $C^2$ in Schikhof's
sense.}\\[2.5mm]
We also provide
alternative characterizations
of $C^{k,\sigma}$-maps
(in the sense of
Bertram et al.)
on open subsets of metrizable spaces,
for $\sigma\in \; ]0,1]$.
Theorem~\ref{dfc} establishes the following characterization.
It is our technically most difficult result,
and its proof relies heavily on a tool
of convenient differential calculus~\cite{KaM},
which has been adapted to non-archimedian analysis
in~\cite{GaLa}.\\[2.5mm]
{\bf Theorem~D.}
\emph{If $\K$ is $\R$ or an ultrametric field
and $E$ is metrizable,
then $f$ is $C^{k,\sigma}$
if and only if $f\circ \gamma\colon \K^{k+1}\to F$
is $C^{k,\sigma}$, for each smooth map
$\gamma\colon \K^{k+1}\to U$.}\\[2.5mm]
Note that neither~$E$ nor~$F$ need to be locally convex
here.
An analogous characterization of
$C^k$-maps was given earlier in \cite[Theorem~12.4]{BGN}.
As a consequence of Theorem~D, the simplified description of
$C^{k,\sigma}$-maps on finite-dimensional
spaces via partial difference quotients
can also be used to
deal with H\"{o}lder differentiable maps
on metrizable spaces.
This may be useful on the way towards
ultrametric (and non-locally convex)
analogues of Boman's
Theorem
(cf.\ \cite[Theorem~2]{Bom}
and \cite[Theorem~12.8]{KaM}),
which characterizes $C^{k,\sigma}$-maps
on open subsets of finite-dimensional (or metrizable)
real locally convex spaces
as those maps which are $C^{k,\sigma}$
along smooth curves.
While the preceding result provided a reduction
to finite-dimensional domains, our next result
(Theorem~\ref{holvsweak})
reduces to the case of a one-dimensional range.\\[2.5mm]
{\bf Theorem~E.}
\emph{If $\K\not=\C$ is locally compact,
$E$ is metrizable and $F$ is locally convex
and Mackey complete,
then $f$ is $C^{k,\sigma}$
if and only if $f$ is weakly $C^{k,\sigma}$,
i.e., $\lambda\circ f\colon U\to \K$
is $C^{k,\sigma}$, for each continuous
linear functional $\lambda\colon F\to\K$.}\\[2.5mm]
We remark that yet another approach
to $C^k$-maps of several variables
has been proposed by De Smedt in~\cite{DS2}.
The $C^1$-maps in the sense of~\cite{DS2}
coincide with the strictly differentiable maps
defined in~\cite{IM2}. If $\K$ is locally
compact, then such maps (on open domains
in finite-dimensional spaces)
coincide with $C^1$-maps in
the sense of Bertram et al.\
(see \cite[Lemma~3.11]{IM2}).\\[2.5mm]
Both the approach to $C^k$-maps
of Schikhof and De Smedt
and the approach of Bertram et al.\
lead to natural topologies
on the function space $C^k(U,F)$,
for~$U$ an open subset of $E=\K^d$.
In Appendix~\ref{apptop},
we show that the two topologies
coincide.
In Appendix~\ref{appreal},
we consider maps from an open set to a real
locally convex space. We show that the $C^{k,\sigma}$-property
can be characterized in terms of the existence
and H\"{o}lder continuity of higher differentials
in this case (Theorem~\ref{viadiff}).\\[2.5mm]
The present studies are part of a larger
project, the goal of which is to transfer
the main ideas of infinite-dimensional
real differential calculus and non-linear
functional analysis into
non-archimedian analysis (and analysis
over arbitrary topological fields).
A survey of the results obtained
so far, with applications
to Lie groups and dynamical systems,
can be found in~\cite{SUR}.
\section{Main concepts, terminology and notation}\label{secnot}
%
%
%
In this section, we compile terminology and notation
concerning differential calculus over topological
fields, together with basic facts.
Most of these facts are easy to take
on faith, and we recommend to skip
the proofs on a first reading.
If desired, the proofs can be looked up in
Appendix~\ref{appbore}.\\[2.5mm]
All topological fields occurring in this article
are assumed Hausdorff and non-discrete;
all topological vector spaces are assumed Hausdorff.
Given a field~$\K$, as usual
we write $\K^\times:=\K\setminus \{0\}$
for its group of invertible elements.
A~valued field is a field~$\K$,
equipped with an absolute value
$|.|\colon \K\to [0,\infty[$
which defines a non-discrete topology on~$\K$.
If $|.|$ satisfies the ultrametric
inequality, we call $(\K,|.|)$
an \emph{ultrametric field}.
Totally disconnected, locally compact
topological fields will be referred to as \emph{local
fields}. It is well known
that each locally compact field
admits an absolute value defining its
topology. We fix such an absolute value,
and thus consider~$\K$ as a valued field.
On $\R$ and $\C$, we shall always use the
usual absolute value.
We write $\N=\{1,2,\ldots\}$
and $\N_0:=\N\cup\{0\}$.
\subsection*{{\normalsize{\boldmath$C^k$}-maps in the sense of
Bertram, Gl\"{o}ckner and Neeb}}
We
recall the approach to $C^k$-maps
between open subsets of topological vector spaces
over a topological field
developed in~\cite{BGN}
(and its extension to maps on
non-open domains from~\cite{IM2}).
More information concerning
this approach can be found in the
survey~\cite{SUR}.
Cf.\ \cite{Ber}
for applications of the corresponding
differential calculus over topological rings
in differential geometry.
We are mostly interested in mappings on open domains,
but some results will hold more generally.\\[3mm]
Let $E$ and $F$ be topological vector
spaces over a topological field~$\K$
and $f\colon  U\to F$ be a map,
defined on a subset $U\sub E$
with dense interior.
Then the directional difference quotient
\[
f^{]1[}(x,y,t)\; :=\;
\frac{f(x+ty)-f(x)}{t}
\]
makes sense for all $(x,y,t)$ in the subset
\[
U^{]1[}\; :=\; \{(x,y,t)\in U\times E\times \K^\times\colon x+ty\in U\}
\]
of $E\times E\times \K$.
To define directional derivatives,
we need to allow also
the value $t=0$.
Hence, we consider
\[
U^{[1]}\;:=\; \{(x,y,t)\in U\times E\times \K
\colon x+ty\in U\}\,.
\]
Then $U^{[1]}=U^{]1[}\cup (U\times E\times \{0\})$, as a disjoint union.
If $U$ is open, then $U^{[1]}$ is an open
subset of the topological $\K$-vector space
$E^{[1]}=E\times E\times \K$.
In the general case, $U^{[1]}\sub E^{[1]}$
has dense interior.
Recursively, we define
$U^{[k]}:=(U^{[1]})^{[k-1]}$
and $U^{]k[}:=(U^{]1[})^{]k-1[}$
for $2\leq k\in \N$.
Then $U^{]k[}$ is dense in
$U^{[k]}$ (see \cite[Remark~1.6]{IM2}).
%
%
\begin{defn}\label{def1}
The map $f\colon U\to F$ is called $C^1_{BGN}$
if $f$ is continuous (i.e., $C^0$, or $C^0_{BGN}$),
and there exists a continuous map
$f^{[1]}\colon U^{[1]}\to F$
which extends
$f^{]1[}\colon U^{]1[}\to F$.
Given $k\in \N$ with $k\geq 2$,
we say that $f$ is
$C^k_{BGN}$ if $f$ is $C^1_{BGN}$ and $f^{[1]}\colon U^{[1]}\to F$
is $C^{k-1}_{BGN}$.
We define
$f^{[k]}:=(f^{[1]})^{[k-1]}\colon
U^{[k]}\to F$ in this case.
The map~$f$ is $C^\infty_{BGN}$
if it is $C^k_{BGN}$ for all $k\in \N_0$.
\end{defn}
Since $U^{]1[}$ is dense in~$U^{[1]}$,
$f^{[1]}$ is unique if
it exists (and likewise each $f^{[k]}$).
%
%
%
\begin{numba}\label{exalin}
For example, every continuous linear
map
$\lambda\colon  E\to F$ is $C^\infty_{BGN}$
with $\lambda^{[1]}(x,y,t)=\lambda(y)$
for all $(x,y,t)\in E\times E\times \K$
(thus also $\lambda^{[1]}$ is continuous linear).
Also each continuous multilinear map
is $C^\infty_{BGN}$ (see \cite{BGN}).
\end{numba}
%
%
%
\begin{numba}\label{chainr}
(Chain Rule).
If $E$, $F$ and $H$ are topological
$\K$-vector spaces, $U\sub E$ and $V\sub F$ are
subsets with dense interior, and $f\colon  U\to V\sub F$,
$g\colon  V\to H$ are $C^k_{BGN}$-maps,
then also the composition
$g\circ f\colon  U\to H$,
$x\mto g(f(x))$ is~$C^k_{BGN}$.
If $k\geq 1$, we have
$(\wh{T}f)(x,y,t):=(f(x),f^{[1]}(x,y,t),t)\in V^{[1]}$
for all $(x,y,t)\in U^{[1]}$, and
%
\begin{equation}\label{formchain}
(g\circ f)^{[1]}(x,y,t)=g^{[1]}(f(x),f^{[1]}(x,y,t),t)\, .
\end{equation}
Thus $(g\circ f)^{[1]}=g^{[1]}\circ \wh{T}f$
with $\wh{T}f\colon U^{[1]}\to V^{[1]}$
(see \cite[Proposition~3.1 and Proposition~4.5]{BGN},
also \cite[\S\,1]{IM2}).
\end{numba}
We recall from \cite[Lemma~4.9]{BGN}
and \cite[\S\,1]{IM2} that
being $C^k$ is a local property.
%
%
\begin{la}\label{Crlocal}
Let $E$ and $F$ be topological $\K$-vector spaces,
and $f\colon  U\to F$ be a map, defined
on a subset $U\sub E$
with dense interior.
Let $k\in \N_0\cup\{\infty\}$. 
If there is an open cover $(U_i)_{i\in I}$
of~$U$ such that
$f|_{U_i}\colon  U_i\to F$ is
$C^k_{BGN}$ for each $i\in I$,
then $f$ is~$C^k_{BGN}$.\Punkt
\end{la}
\subsection*{{\normalsize{\boldmath $C^k$}-maps in the sense
of Schikhof and De Smedt}}
In this section, we give a definition
of $C^k$-maps of several variables
based on continuous extensions
to certain partial difference quotient
maps,
which generalizes special cases
considered by Schikhof \cite[\S\,84]{Sch}
and De Smedt \cite{DSm}.
Our notation differs from the one used
in~\cite{DSm} and~\cite{Sch},
because we find it more convenient
to use multi-indices in
higher dimensions.
\begin{numba}
Until Remark~\ref{onlynois},
let $\K$ be a topological field,
$d\in \N$, $U\sub \K^d$
be an open subset
(where $\K^d$ is equipped with the product
topology), and $F$ be a topological
$\K$-vector space.
As usual, for $i\in \{1,\ldots, d\}$
we set $e_i:=(0,\ldots,0,1,0,\ldots, 0)\in \K^d$,
with $i$-th entry~$1$.
\end{numba}
%
%
\begin{numba}\label{notconv}
As usual, given a ``multi-index'' $\alpha=(\alpha_1,\ldots,\alpha_d)
\in \N_0^d$,
we write $|\alpha|:=\sum_{i=1}^d \alpha_i$.
The definition of a $C^k$-map $f\colon U\to F$ in the
sense of Schikhof and De Smedt
will involve a certain continuous extension
$f^{<\alpha>}$ of a partial difference
quotient map $f^{>\alpha<}$ corresponding to
each multi-index $\alpha\in \N_0^d$
such that $|\alpha|\leq k$.
It is convenient to define
the domains $U^{<\alpha>}$ and $U^{>\alpha<}$
of these mappings first.
They will be subsets
of $\K^{d+|\alpha|}$.
It is useful to write elements
$x\in \K^{d+|\alpha|}$ in the form
$x=(x^{(1)}, x^{(2)},\ldots, x^{(d)})$,
where $x^{(i)}\in \K^{1+\alpha_i}$
for $i\in \{1,\ldots, d\}$.
We write $x^{(i)}=
(x^{(i)}_0,x^{(i)}_1,\ldots, x^{(i)}_{\alpha_i})$
with $x^{(i)}_j\in \K$
for $j\in \{0,\ldots,\alpha_i\}$.
\end{numba}
\begin{numba}
Given $\alpha\in \N_0^d$, we now define
$U^{<\alpha>}$ as the set of all
$x\in \K^{d+|\alpha|}$
such that, for all
$i_1\in \{0,1,\ldots,\alpha_1\},\ldots,
i_d \in \{0, 1,\ldots,\alpha_d\}$,
we have
\[
(x^{(1)}_{i_1},\ldots,
x^{(d)}_{i_d})\in U\, .
\]
We let $U^{>\alpha<}$ be the set of all
$x\in U^{<\alpha>}$ such that,
for all $i\in \{1,\ldots, d\}$
and $0\leq j<k\leq \alpha_i$,
we have $x^{(i)}_j\not=x^{(i)}_k$.
It is easy to see that~$U^{<\alpha>}$ and~$U^{>\alpha<}$
are open in $\K^{d+|\alpha|}$
and~$U^{>\alpha<}$
is dense in~$U^{<\alpha>}$.
\end{numba}
%
%
\begin{example}\label{mimpex}
If $U=U_1\times\cdots\times U_d$
with open sets $U_i\sub \K$,
then simply
%
\begin{equation}\label{spcas}
U^{<\alpha>}\; =\; U_1^{1+\alpha_1}\times
U_2^{1+\alpha_2}\times\cdots
\times U_d^{1+\alpha_d}\,.
\end{equation}
Only this case (in fact only special
cases thereof) was considered
in~\cite{DSm}~and~\cite{Sch}.
\end{example}
\begin{rem}
A simple induction on $|\alpha|$ shows that
the sets $U^{<\alpha>}$ can be defined
alternatively by recursion on $|\alpha|$,
as follows: Set $U^{<0>}:=U$.
Given $\alpha\in \N_0^d$ such that
$|\alpha|\geq 1$, pick $\beta\in \N_0^d$
such that $\alpha=\beta+e_i$
for some $i\in \{1,\ldots, d\}$.
Then
$U^{<\alpha>}$
is the set of all elements $x\in \K^{d+|\alpha|}$
such that\linebreak
$(x^{(1)},\ldots, x^{(i-1)},
x^{(i)}_0,x^{(i)}_1,\ldots,x^{(i)}_{\alpha_i-1},
x^{(i+1)},\ldots, x^{(d)})\in U^{<\beta>}$
holds as well as\linebreak
$(x^{(1)},\ldots, x^{(i-1)},
x^{(i)}_{\alpha_i},x^{(i)}_1,\ldots,x^{(i)}_{\alpha_i-1},
x^{(i+1)},\ldots, x^{(d)})\in U^{<\beta>}$.
\end{rem}
We now define certain mappings
$f^{>\alpha<}\colon U^{>\alpha<}\to F$ and show afterwards
that they can be interpreted
as partial difference
quotient maps.
\begin{defn}
We set $f^{>0<}:=f$.
Given a multi-index $\alpha\in \N_0^d$
such that $|\alpha|\geq 1$,
we define $f^{>\alpha<}(x)$ as the sum
%
\begin{equation}\label{badbd}
\sum_{j_1=0}^{\alpha_1}\cdots
\sum_{j_d=0}^{\alpha_d}
\left(\prod_{k_1\not=j_1} \frac{1}{x^{(1)}_{j_1}-x^{(1)}_{k_1}}
\cdot\ldots\cdot
\prod_{k_d\not=j_d}\frac{1}{x^{(d)}_{j_d}-x^{(d)}_{k_d}}\right)
f(x^{(1)}_{j_1},\ldots, x^{(d)}_{j_d})
\end{equation}
for $x\in U^{>\alpha<}$,
using the notational conventions
from {\bf\ref{notconv}}.
The products are taken over all
$k_\ell\in \{0,\ldots, \alpha_\ell\}$
such that $k_\ell\not=j_\ell$,
for $\ell\in \{1,\ldots, d\}$.
\end{defn}
The map $f^{>\alpha<}$ has
important symmetry properties.
%
%
\begin{la}\label{presymm}
Assume that $\alpha\in \N_0^d$,
$i\in \{1,\ldots, d\}$
and $\pi$ is a permutation of\linebreak
$\{0,1,\ldots, \alpha_i\}$.
Then
$(x^{(1)},\ldots, x^{(i-1)},
x^{(i)}_{\pi(0)},\ldots, x^{(i)}_{\pi(\alpha_i)},
x^{(i+1)},\ldots, x^{(d)})\in U^{>\alpha<}$
for each $x\in U^{>\alpha<}$,
and
%
\begin{equation}\label{preeqsmm}
f^{>\alpha<}(x^{(1)},\ldots, x^{(i-1)},
x^{(i)}_{\pi(0)},\ldots, x^{(i)}_{\pi(\alpha_i)},
x^{(i+1)},\ldots, x^{(d)})\;=\; f^{>\alpha<}(x)\, .
\end{equation}
\end{la}
The following lemma
shows that $f^{>\alpha<}$ can indeed be
interpreted as a partial difference quotient map.
%
%
\begin{la}\label{interpr}
For each $i\in \{1,\ldots, d\}$ and $x\in U^{>e_i<}$,
the element $f^{>e_i<}(x)$ is given by
\[
\frac{f(x^{(1)},\ldots, x^{(i-1)},x^{(i)}_0, x^{(i+1)},\ldots, x^{(d)})
-
f(x^{(1)},\ldots, x^{(i-1)},x^{(i)}_1,
x^{(i+1)},\ldots, x^{(d)})}{x^{(i)}_0-x^{(i)}_1}.
\]
If $\alpha\in \N_0^d$ such that $|\alpha|\geq 2$,
let
$\beta\in \N_0^d$ be a multi-index such that
$\alpha=\beta+e_i$
for some
$i\in \{1,\ldots,d\}$. Then $f^{>\alpha<}(x)$ is given by
%
\begin{eqnarray}
\hspace*{-5mm}\lefteqn{\frac{1}{x^{(i)}_0-x^{(i)}_{\alpha_i}}\cdot
\Big(f^{>\beta<}(x^{(1)},\ldots, x^{(i-1)},
x^{(i)}_0, x^{(i)}_1, \ldots, x^{(i)}_{\alpha_i-1},
x^{(i+1)},\ldots, x^{(d)})}\quad\quad\quad \notag\\
& & \;\;\; -
f^{>\beta<}(x^{(1)},\ldots, x^{(i-1)},
x^{(i)}_{\alpha_i}, x^{(i)}_1, \ldots, x^{(i)}_{\alpha_i-1},
x^{(i+1)},\ldots, x^{(d)})\Big)
\label{schlmm}
\end{eqnarray}
for all $x\in U^{>\alpha<}$.
\end{la}
%
%
\begin{defn}\label{defSDS}
We say that $f$ is $C^0_{SDS}$ if it is continuous,
and define $f^{<0>}:=f$ in this case.
Recursively, given an integer $k\geq 1$
we say that
$f$ is $C^k_{SDS}$
if $f$ is $C^{k-1}_{SDS}$
and, for each multi-index
$\alpha\in \N_0^d$ such that $|\alpha|=k$,
there exists
a continuous map
$f^{<\alpha>}\colon U^{<\alpha>}\to F$
such that $f^{<\alpha>}|_{U^{>\alpha<}}=f^{>\alpha<}$.
As usual, $f$ is called $C^\infty_{SDS}$
if $f$ is $C^k_{SDS}$ for each $k\in \N_0$.
\end{defn}
Since $U^{>\alpha<}$ is dense in
$U^{<\alpha>}$, the continuous extension
$f^{<\alpha>}$ of~$f^{>\alpha<}$
is unique whenever it exists.
We readily deduce from Lemma~\ref{presymm}:
%
%
\begin{la}\label{partsymm}
Let $f$ be a $C^k$-map
for some $k\in \N$, $\alpha\in \N_0^d$
with $|\alpha|=k$,\linebreak
$i\in \{1,\ldots, d\}$
and $\pi$ be a permutation of
$\{0,1,\ldots, \alpha_i\}$.
Then
%
\begin{equation}\label{nweq}
(x^{(1)},\ldots, x^{(i-1)},
x^{(i)}_{\pi(0)},\ldots, x^{(i)}_{\pi(\alpha_i)},
x^{(i+1)},\ldots, x^{(d)})\; \in \; U^{<\alpha>}
\end{equation}
for each $x\in U^{<\alpha>}$, and
%
\begin{equation}\label{eqsmm}
f^{<\alpha>}(x^{(1)},\ldots, x^{(i-1)},
x^{(i)}_{\pi(0)},\ldots, x^{(i)}_{\pi(\alpha_i)},
x^{(i+1)},\ldots, x^{(d)})\;=\; f^{<\alpha>}(x)\, .
\end{equation}
\end{la}
The following variant of Lemma~\ref{interpr}
is available for~$f^{<\alpha>}$.
%
%
\begin{la}\label{extmore}
Let $f$ be a $C^k_{SDS}$-map
for an integer $k\geq 2$,
$\alpha\in \N_0^d$ such that $|\alpha|=k$,
and $\beta\in \N_0^d$ such that
$\alpha=\beta+e_i$ for some
$i\in \{1,\ldots, d\}$.
Then $f^{<\alpha>}(x)$ is given by
%
%
\begin{eqnarray}
\hspace*{-5mm}
\lefteqn{\frac{1}{x^{(i)}_0-x^{(i)}_{\alpha_i}}
\cdot \Big(f^{<\beta>}(x^{(1)},\ldots, x^{(i-1)},
x^{(i)}_0, x^{(i)}_1, \ldots, x^{(i)}_{\alpha_i-1},
x^{(i+1)},\ldots, x^{(d)})}\qquad\quad \notag\\
& & \;\;\; -
f^{<\beta>}(x^{(1)},\ldots, x^{(i-1)},
x^{(i)}_{\alpha_i}, x^{(i)}_1, \ldots, x^{(i)}_{\alpha_i-1},
x^{(i+1)},\ldots, x^{(d)})\Big)\label{schlmm2}
\end{eqnarray}
for all $x\in U^{<\alpha>}$
such that $x^{(i)}_0\not=x^{(i)}_{\alpha_i}$.
\end{la}
%
%
%
\begin{rem}\label{onlynois}
If $U\sub \K^d$ is a subset
(possibly with empty interior)
of the form $U=U_1\times\cdots\times U_d$,
where $U_i\sub \K$ is a non-empty subset
without isolated points for $i\in \{1,\ldots, d\}$,
then Definition~\ref{defSDS} can be used
just as well to define $C^k_{SDS}$-maps
$f\colon U\to F$.
\end{rem}
If $d=1$, we also write
$f^{<j>}$ in place of
$f^{<j e_1>}$, as in~\cite[\S\,6]{BGN}.
\subsection*{{\normalsize Seminorms and gauges}}
Gauges on topological vector spaces over
valued fields were introduced in~\cite{IM2} as a substitute
for continuous seminorms when dealing with a
general topological vector space,
the topology of which need not come from a family
of continuous seminorms
(cf.\ \cite[\S\,6.3]{Jar} for the real case).
We only recall some essentials here; see~\cite{IM2}
for further information.
\begin{defn}
Let $E$ be a topological
vector space over a valued field
$(\K,|.|)$. A \emph{gauge on~$E$}
is a map $q \colon E\to [0,\infty[$
(also written $\|.\|_q := q$)
which satisfies
$q (t x) = |t|q(x)$ for all $t\in \K$
and $x\in E$,
and such that $B_r^q(0)$
is a $0$-neighbourhood
for each $r>0$,
where
$B_r^q(x):=\{y\in E\colon \|y-x\|_q<r\}$
for all $x\in E$ and $r>0$.
We also define
$\wb{B}_r^q(x):=\{y\in E\colon \|y-x\|_q\leq r\}$.
If $(E,\|.\|)$ is a normed
space, we relax notation and write
$B^E_r(x):=B^{\|.\|}_r(x)$.
\end{defn}
In~\cite{IM2},
only upper semicontinuous gauges~$q$
were considered, i.e., it was required that
$B^q_r(0)$ is an open $0$-neighbourhood,
for each $r>0$.
%
%
\begin{rem}\label{minkow}
Typical examples
of gauges are Minkowski functionals~$\mu_U$
of balanced, open $0$-neighbourhoods~$U$
in a topological vector space~$E$
over a valued field~$\K$;
these are upper semicontinuous
(see \cite[Remark\,1.21]{IM2}).
Here $U\sub E$ is called
\emph{balanced} if $tU\sub U$
for all $t\in \K$ such that $|t|\leq 1$.
The Minkowski functional
is
$\mu_U\colon E\to [0,\infty[$,
$x\mto \inf\{|t|\colon \mbox{$t\in \K^\times$
with $x\in tU$}\}$.
\end{rem}
%
%
\begin{rem}\label{remfk}
Note that gauges
need not satisfy the triangle inequality.
But we still have a certain substitute:
Given a gauge $q\colon E\to [0,\infty[$,
there always exists a gauge
$p\colon E\to [0,\infty[$
such that
%
%
\begin{equation}\label{faketri}
\|x+y\|_q\;\leq\; \|x\|_p+\|y\|_p
\quad\mbox{for all $x,y\in E$}
\end{equation}
(cf.\ \cite[Lemma~1.29]{IM2}). We shall refer to
(\ref{faketri}) as the \emph{fake
triangle inequality}.
\end{rem}
As in the case of continuous seminorms, it
frequently suffices to consider
a sufficiently large set of gauges:
%
%
\begin{defn}\label{funsys}
A set $\Gamma$ of gauges
on a topological $\K$-vector space~$E$
is called a \emph{fundamental system
of gauges} if each $0$-neighbourhood in~$E$
contains some finite intersection
of balls of the form $B_r^q(0)$,
with $q\in \Gamma$ and $r>0$.
\end{defn}
Cf.\ \cite[Lemma 1.24]{IM2}
for the following lemma, which is useful
to determine fundamental systems
of gauges.
%
%
\begin{la}\label{findfund}
Let $p,q\colon E\to [0,\infty[$
be gauges on a topological vector space~$E$ over
a valued field $\K$.
If there exist $r,s>0$ such that
$B^q_s(0)\sub B^p_r(0)$,
then
\[
p \; \leq rs^{-1}|a|^{-1}q\, ,
\]
for each $a\in \K^\times$ such that $|a|<1$.
In particular, $p\leq C q$ for some $C>0$.\,\Punkt
\end{la}
%
%
\begin{rem}\label{examfun}
Combining Remark~\ref{minkow}
and Lemma~\ref{findfund}, it is easy to see
that upper semicontinuous gauges form a
fundamental system of gauges,
for each topological vector space over
a valued field (cf.\ also \cite[Remark~1.21]{IM2}).
In the real case, continuous
gauges form a fundamental system
(cf.\ \cite[\S\,6.4]{Jar}).
\end{rem}
\begin{examples}
Given $r\in \;]0,1]$,
a gauge $q\colon E\to [0,\infty[$
is called an $r$-seminorm if
$q(x+y)^r\leq q(x)^r+q(y)^r$
for all $x,y\in E$.
If, furthermore,
$q(x)=0$ if and only if $x=0$,
then~$q$ is called an $r$-norm
(cf.\ \cite[\S\,6.3]{Jar}
for the real case).
For examples of $r$-normed spaces over~$\R$
and more general non-locally convex
real topological vector spaces,
the reader is referred to \cite[\S\,6.10]{Jar}
and~\cite{Kal}.
For $\K$ a valued field,
the simplest
examples are the spaces
$\ell^p(\K)$
of all $x=(x_n)_{n\in\N}
\in \K^\N$
such that $\|x\|_p:=\sqrt[p]{\sum_{n=1}^\infty|x_n|^p}<\infty$,
for $p\in \;]0,1[$.
Then $\|.\|_p$
is a $p$-norm
on $\ell^p(\K)$ defining a Hausdorff
vector topology on this space
(and thus $\{\|.\|_p\}$ is a fundamental
system of gauges).
\end{examples}
\subsection*{{\normalsize Bounded sets and bounded maps}}
Let $E$ be a topological
vector space over a topological
field~$\K$.
Recall that a subset $B\sub E$
is called \emph{bounded}
if, for each $0$-neighbourhood
$U\sub E$, there exists a $0$-neighbourhood
$V\sub \K$ such that $VB\sub U$.
If $\K$ is a valued field,
we can test boundedness using gauges.
%
%
\begin{la}\label{bdviagau}
Let $E$ be a topological vector space
over a valued field~$\K$.
Then a subset $B\sub E$ is bounded
if and only if the set $q(B)\sub \R$ is bounded,
for each gauge $q$ on~$E$.
\end{la}
It suffices to show $\,\sup \|B\|_q<\infty$
for $q$ in a fundamental system
of gauges. In Section~\ref{secmetr},
$C^k$-maps with bounded difference
quotients will play a vital~role.
\begin{defn}
Let $E$ be a topological vector space
over a topological field~$\K$.
\begin{itemize}
\item[(a)]
If~$X$ is a topological space,
then $BC(X,E)$ denotes the set
of all continuous maps $\gamma\colon X\to E$
whose image $\gamma(X)$ is bounded in~$E$.
Then $BC(X,E)$ is a vector subspace
of $E^X$.
\item[(b)]
If $k\in \N_0\cup\{\infty\}$ and
$U\sub \K$ is a subset without isolated points,
we let
$BC^k(U,E)$ be the space of all $C^k_{SDS}$-maps
$\gamma\colon U\to E$ such that
$\gamma^{<j>}\in BC(U^{j+1},E)$
for all $j\in \N_0$ such that $j\leq k$.
\end{itemize}
\end{defn}
We mention that $BC^k(U,E)$ can be made a topological
vector space~\cite[Definition~1.2]{GaLa},
but we shall not use this topology.
If~$\K$ is a valued field,
$j\in \N_0$ with $j\leq k$ and
$q$ a gauge on~$E$, we define
%
%
\begin{equation}\label{nogau}
\|\gamma^{<j>}\|_{q,\infty}\; :=\;
\sup\{\|\gamma^{<j>}(x)\|_q\colon x\in U^{j+1}\}\quad
\mbox{for $\gamma\in BC^k(U,E)$.}
\end{equation}
\subsection*{{\normalsize H\"{o}lder continuity}}
Using gauges, we now define a version of H\"{o}lder
continuous maps and record their basic properties.
%
\begin{defn}\label{deflipcts}
Let $E$ and $F$
be topological vector spaces over a valued field~$\K$,
and $U\sub E$ be a subset.
A map $f\colon U\to F$
is called \emph{H\"{o}lder continuous
of exponent $\sigma\in\;]0,\infty[$}
(or $C^{0,\sigma}$)
if, for every $x_0\in U$ and gauge $q$
on $F$, there exists a gauge~$p$ on~$E$
and a neighbourhood $V\sub U$ of~$x_0$
such that
%
%
\begin{equation}\label{hoelcon}
\|f(y)-f(x)\|_q \; \leq \; \big(\|y-x\|_p\big)^\sigma
\quad
\mbox{for all $\, x,y\in V$.}
\end{equation}
$C^{0,1}$-maps are also called
\emph{Lipschitz continuous.}
\end{defn}
We remind the reader that
H\"{o}lder exponents $\sigma>1$
are meaningful in non-archimedian analysis
(see \cite[Exercise~26.B]{Sch} for an
instructive example).
We are mostly interested in H\"{o}lder
exponents $\sigma\in \;]0,1]$,
but some of the results are valid
just as well for $\sigma>1$.
%
%
%
\begin{la}\label{baseHglob}
Let $E$, $F$ and $H$ be topological
vector spaces over a valued field~$\K$,
$U\sub E$ and $V\sub F$ be subsets,
$f\colon U\to V\sub F$
and $g\colon V\to H$ be maps,
and $\sigma,\tau > 0 $.
Then the following holds:
\begin{itemize}
\item[\rm (a)]
If $f$ is $C^{0,\sigma}$,
then $f$ is continuous.
\item[\rm (b)]
If $f$ is $C^{0,\sigma}$ and $\sigma\geq \tau$,
then $f$ is also $C^{0,\tau}$.
\item[\rm (c)]
If $f$ is $C^{0,\sigma}$ and $g$ is $C^{0,\tau}$,
then $g\circ f$ is $C^{0,\sigma\cdot\tau}$.
\item[\rm (d)]
If $U$ has dense interior
and $f$ is $C^1_{BGN}$, then $f$ is Lipschitz continuous.
\end{itemize}
\end{la}
In connection with Part\,(d) of the preceding lemma,
note that $\id\colon \K\to\K$, $x\mto x$
is $C^\infty_{BGN}$ but not $C^{0,\sigma}$ for any $\sigma>1$.\\[2.5mm]
If $f$ is not H\"{o}lder continuous,
then pairs of points with pathological
behaviour
can always be chosen in a given dense set.
This will become essential later.
%
%
\begin{la}\label{densehoel}
Let $E$ and $F$ be topological
vector spaces over a valued field~$\K$,
$f\colon U\to F$ be a continuous mapping
on a subset
$U\sub E$,
$D\sub U$ be a dense subset,
and $\sigma > 0$.
If $f$ is not $C^{0,\sigma}$,
then there exists $x_0\in U$
and a gauge~$q$ on~$F$ such that,
for each neighbourhood $V\sub U$
of~$x_0$ and gauge~$p$ on~$E$,
there exist $x,y\in V\cap D$ such that
$\|f(x)-f(y)\|_q>(\|x-y\|_p)^\sigma$.
\end{la}
\subsection*{{\normalsize Two approaches to
H\"{o}lder differentiable maps}}
We define $k$ times H\"{o}lder differentiable maps
and record some properties.
\begin{defn}
Let $\K$ be a valued field,
$E$ and $F$ be topological $\K$-vector spaces,
$U\sub E$ be a subset
with dense interior,
and $f\colon U\to F$ be a mapping.
Let $k\in \N_0\cup\{\infty\}$
and $\sigma>0$.
We say that $f$ is $C^{\,k,\sigma}_{BGN}$
if $f$ is $C^k_{BGN}$
and $f^{[j]}\colon U^{[j]}\to F$
is $C^{0,\sigma}$
for all $j\in \N_0$
such that $j\leq k$
(where $f^{[0]}:=f$).
\end{defn}
%
%
\begin{rem}\label{autolip}
Note that, if $\sigma\in \;]0,1]$,
then $f^{[j]}$ is $C^1_{BGN}$ for $j<k$
and hence automatically~$C^{0,\sigma}$,
by Lemma~\ref{baseHglob}~(b) and~(d).
In particular, if $\sigma\in \;]0,1]$,
then $f$ is $C^{\infty,\sigma}_{BGN}$ if and only if
$f$ is $C^\infty_{BGN}$.
And if $k$ is finite, then
only $f^{[k]}$ requires attention.
\end{rem}
\begin{rem}
$k$ times Lipschitz differentiable mappings
($C^{k,1}$-maps)
form a particularly nice class
of maps (see~\cite{IM2}).
Notably, \cite[Theorem~5.2]{IM2} provides
an implicit function theorem
for $C^{k,1}$-maps
from arbitrary topological vector spaces
to Banach spaces, for each valued field
and $k\in \N\cup\{\infty\}$.
\end{rem}
We need some basic information on mappings
into direct products.
\begin{la}\label{prodhoel}
Let $E$ be a topological $\K$-vector space
over a valued field~$\K$, $(F_i)_{i\in I}$
be a family of topological $\K$-vector spaces
and $f
\colon U\to F$ be a map into
$F:=\prod_{i\in I}F_i$,
defined on a non-empty subset $U\sub E$
with dense interior. Let $k\in \N_0$ and $\sigma>0$;
for $i\in I$, let $\pr_i\colon F\to F_i$ be the projection.
Then $f$ is $C^{k,\sigma}_{BGN}$ if and only if
each of its components
$f_i:=\pr_i\circ f \colon U\to F_i$ is $C^{k,\sigma}_{BGN}$.
\end{la}
The Chain Rule is available in the following form.
%
%
\begin{la}\label{chainLip}
Let $\K$ be a valued field,
$E$, $F$ and $H$ be topological $\K$-vector spaces,
$U\sub E$ and $V\sub F$ be subsets with dense interior,
$\sigma\in \;]0,1]$, $\tau >0$,
$f\colon U\to V$ be $C^{\,k,\sigma}_{BGN}$,
and $g\colon V\to H$ be
$C^{\,k,\tau}_{BGN}$.
Then $g\circ f\colon U\to H$
is $C^{\,k,\sigma\cdot\tau}_{BGN}$.
\end{la}
The following variant even holds if $\sigma>1$:
%
\begin{la}\label{chainlin}
Let $\K$ be a valued field,
$E$, $F$ and $H$ be topological $\K$-vector spaces,
$\lambda\colon F\to H$ be continuous linear,
$U\sub E$ be a subset with dense interior,
$\sigma> 0$
and $f\colon U\to F$ be a
$C^{\,k,\sigma}_{BGN}$-map, where $k\in \N_0$.
Then $\lambda \circ f\colon U\to H$
is $C^{\,k,\sigma}_{BGN}$,
and $(\lambda\circ f)^{[k]}=\lambda
\circ f^{[k]}$.
\end{la}
%
%
%
\begin{la}\label{hoelloc}
Let $E$ and $F$ be topological vector spaces over a
valued field~$\K$,
and $f\colon  U\to F$ be a map, defined
on a subset $U\sub E$
with dense interior.
Let $k\in \N_0\cup\{\infty\}$
and $\sigma\in \;]0,1]$.
If there is an open cover $(U_i)_{i\in I}$
of~$U$ such that
$f|_{U_i}\colon  U_i\to F$ is
$C^{k,\sigma}_{BGN}$ for each $i\in I$,
then $f$ is~$C^{k,\sigma}_{BGN}$.
\end{la}
H\"{o}lder differentiable maps
can also be defined using the approach of Schikhof and De Smedt.
\begin{defn}
Let $\K$ be a valued field, $d\in \N_0$,
$F$ be a topological $\K$-vector space
and $f\colon U\to F$ be a mapping,
where $U\sub \K^d$ is open or
$U=U_1\times \cdots \times U_d$
for certain sets $U_1,\ldots, U_d \sub \K$ without
isolated points.
Let $k\in \N_0\cup\{\infty\}$
and $\sigma > 0 $.
We say that $f$ is $C^{\,k,\sigma}_{SDS}$
if $f$ is $C^k_{SDS}$
and $f^{<\alpha>}\colon U^{<\alpha>}\to F$
is $C^{0,\sigma}$
for all $\alpha \in \N_0^d$
such that $|\alpha|\leq k$
$($where $f^{<0>}:=f)$.
\end{defn}
\section{{\boldmath$C^k_{BGN}$}-maps and
{\boldmath$C^k_{SDS}$}-maps coincide}\label{seccoinc}
%
%
%
In this section, we show that the approach of
Bertram, Gl\"{o}ckner and Neeb
and the approach of Schikhof and De Smedt
give rise to the same classes of~$C^k$-maps
and $C^{k,\sigma}$-maps
on open domains (and more generally).
Throughout the section,
$\K$ is a topological field,
$d\in \N$ and
$f\colon U\to F$ a map
to a topological $\K$-vector space~$F$,
where $U\sub \K^d$ is open or of the
form $U=U_1\times\cdots\times U_d$
for certain sets $U_1,\ldots, U_d\sub \K$
with dense interior.
%
%
\begin{thm}\label{mainthm}
The following holds
for each $k\in \N_0\cup\{\infty\}$:
\begin{itemize}
\item[\rm (a)]
$f$ is
$C^k_{SDS}$
if and only if $f$ is $C^k_{BGN}$.
\item[\rm(b)]
If $\K$ is a valued field
and $\sigma\in \;]0,1]$,
then $f$ is
$C^{\,k,\sigma}_{SDS}$
if and only if $f$ is $C^{\,k,\sigma}_{BGN}$.
\end{itemize}
\end{thm}
Various lemmas are useful for the proof
of Theorem~\ref{mainthm}.
%
%
\begin{la}\label{enindSDS}
Let $k\in \N_0$.
\begin{itemize}
\item[\rm(a)]
If $f\colon \K^d\supseteq U\to F$
is $C^1_{SDS}$ and $f^{<e_i>}$
is $C^k_{SDS}$ for each $i\in \{1,\ldots, d\}$,
then $f$ is $C^{k+1}_{SDS}$.
\item[\rm (b)]
Let $\K$ be a valued field
and $\sigma > 0$.
If $f\colon \K^d\supseteq U\to F$
is $C^{\,1,\sigma}_{SDS}$ and $f^{<e_i>}$
is $C^{\,k,\sigma}_{SDS}$ for each $i\in \{1,\ldots, d\}$,
then $f$ is $C^{\,k+1,\sigma}_{SDS}$.
\end{itemize}
\end{la}
\begin{proof}
Given $\alpha\in \N_0^d$ such that $1\leq |\alpha|\leq k+1$,
there is $i\in \{1,\ldots, d \}$
such that $\alpha_i>0$. Then $\beta:=\alpha-e_i\in
\N_0^d$ and $|\beta|=|\alpha|-1\leq k$.
Write $\beta=(\beta_1,\ldots, \beta_d)$
and set $\beta':=(\beta_1,\ldots,\beta_i,0,\beta_{i+1},\ldots,
\beta_d)\in \N_0^{d+1}$.
For
$x\in U^{>\alpha<}$,
using the notational conventions from~{\bf\ref{notconv}}, we have
\[
f^{>\alpha<}(x)=
(f^{>e_i<})^{>\beta'\!<}
(x^{(1)};\ldots; x^{(i-1)};
x^{(i)}_0,x^{(i)}_2,\ldots, x^{(i)}_{\alpha_i};
x^{(i)}_1;
x^{(i+1)};\ldots; x^{(d)}),
\]
as is clear from the definitions.
Thus
\[
f^{<\alpha>}(x):=
(f^{<e_i>})^{<\beta'\!>}
(x^{(1)};\ldots; x^{(i-1)};
x^{(i)}_0,x^{(i)}_2,\ldots, x^{(i)}_{\alpha_i};
x^{(i)}_1;
x^{(i+1)};\ldots; x^{(d)})
\]
for $x\in U^{<\alpha>}$
defines a continuous (resp., $C^{\,0,\sigma}$-)
extension
$f^{<\alpha>}\colon U^{<\alpha>}\to F$
of $f^{>\alpha<}$,
whenever $|\alpha|\leq k+1$.
Therefore, $f$ is $C^{k+1}_{SDS}$
(resp., $C^{\,k+1,\sigma}_{SDS}$).
\end{proof}
The next lemma establishes one implication
in Theorem~\ref{mainthm} (a) and (b).
Note that $\sigma$ need not be $\leq 1$
here.
%
%
\begin{la}\label{BGNSDS}
If $f$ is $C^k_{BGN}$
for some $k\in \N_0$
$($resp., $C^{\,k,\sigma}_{BGN}$
if $\K$ is a valued field
and $\sigma>0 )$,
then $f$ is $C^k_{SDS}$
$($resp., $C^{\,k,\sigma}_{SDS})$.
\end{la}
\begin{proof}
The proof is by induction
on $k\in \N_0$. The case $k=0$ is trivial.
If $k\geq 1$,
let $i\in \{1,\ldots, n\}$.
For each $x\in U^{>e_i<}$,
we then have
\[
f^{>e_i<}(x)\;=\; f^{[1]}(x^{(1)},\ldots,x^{(i-1)},
x^{(i)}_1,x^{(i+1)},\ldots, x^{(d)}; e_i; \, x^{(i)}_0-x^{(i)}_1)\,.
\]
Thus
%
\begin{equation}\label{RHint}
f^{<e_i>}(x)\;:=\; f^{[1]}(x^{(1)},\ldots,x^{(i-1)},
x^{(i)}_1,x^{(i+1)},\ldots, x^{(d)}; e_i; \, x^{(i)}_0-x^{(i)}_1)
\end{equation}
for $x\in U^{<e_i>}$
defines a continuous extension
of $f^{>e_i<}$
and hence $f$ is $C^1_{SDS}$,
with $f^{<e_i>}$ as just described.
Note that the right hand side of (\ref{RHint})
expresses $f^{<e_i>}$
as a composition of the $C^{k-1}_{BGN}$-map
(resp., $C^{\,k-1,\sigma}_{BGN}$-map)
$f^{[1]}$ and the restriction
of a map $\K^{d+1}\to \K^d\times\K^d\times\K$
which is
continuous affine-linear and hence $C^\infty_{BGN}$.
By the Chain Rule (resp., Lemma~\ref{chainLip}),
$f^{<e_i>}$ is $C^{k-1}_{BGN}$
(resp., $C^{\,k-1,\sigma}_{BGN}$)
and hence $C^{k-1}_{SDS}$
(resp., $C^{\,k-1,\sigma}_{SDS}$),
by induction.
Hence $f$ is $C^{k+1}_{SDS}$ (resp.,
$C^{\,k+1,\sigma}_{SDS}$), by Lemma~\ref{enindSDS}.
\end{proof}
%
%
%
\begin{la}\label{preSDS}
If $f$ is $C^k_{SDS}$ for some $k\in \N_0$
$($resp., if $\K$ is a valued field
and $f$ is $C^{\,k,\sigma}_{SDS}$
for some $\sigma>0 )$,
then $f^{<\alpha>}$ is $C^{k-|\alpha|}_{SDS}$
$($resp., $C^{k-|\alpha|,\sigma}_{SDS})$
for each $\alpha\in \N_0^d$ such that
$|\alpha|\leq k$.
\end{la}
%
%
\begin{rem}\label{remeazy}
The proof of Lemma~\ref{preSDS}
will provide the following formula
for $(f^{<\alpha>})^{<\beta>}$,
if $\alpha\in \N_0^d$ such that $|\alpha|\leq k$
and $\beta\in \N_0^{d+|\alpha|}$
such that $|\beta|\leq k-|\alpha|$:
%
\begin{equation}\label{simpfml}
(f^{<\alpha>})^{<\beta>}\;=\;f^{<\alpha+\bar{\beta}>}\,,
\end{equation}
where $\bar{\beta}\in \N_0^d$ is defined by
$\bar{\beta}_j:=\sum_{i=s_j}^{s_{j+1}-1}\beta_i$
for $j\in \{1,\ldots, d\}$,
with $s_j:=j+\sum_{i=1}^{j-1}\alpha_i$
and $s_{d+1}:=d+|\alpha|+1$.
\end{rem}
\textbf{Proof of Lemma~\ref{preSDS}.}
Given $\alpha\in \N_0^d$,
we first show by induction on\linebreak
$\ell \in \{0,\ldots,k-|\alpha|\}$ that
%
%
\begin{equation}\label{simpfml2}
(f^{>\alpha<})^{>\beta<}\;=\;f^{>\alpha+\bar{\beta}<}\,,
\end{equation}
for each $\beta\in \N_0^{d+|\alpha|}$
such that $|\beta|=\ell$,
where~$\bar{\beta}$ (and $s_1,\ldots, s_{d+1}$)
are as in Remark~\ref{remeazy}.
The case $\ell=0$ being trivial,
let us assume now that (\ref{simpfml2})
holds for some $\beta$ with $|\beta|<k-|\alpha|$.
For each $i \in \{1,\ldots, d+|\alpha|\}$,
we have to show that~(\ref{simpfml2})
holds with $\beta$ replaced by $\gamma:=\beta+e_i$
and $\bar{\beta}$ replaced by the corresponding~$\bar{\gamma}$.
There is a unique
$j \in \{1,\ldots, d\}$
such that $s_j\leq i <s_{j+1}$.
Then $\bar{\gamma}=\bar{\beta}+e_j$.
Given $x\in (U^{>\alpha<})^{>\gamma<}
=(U^{>\alpha<})^{>\beta+e_i<}\sub (\K^{d+|\alpha|})^{>\beta+e_i<}$,
we abbreviate
$y:=(x^{(1)},\ldots, x^{(s_j-1)})$
and $z:=(x^{(s_j+1)},\ldots, x^{(d+|\alpha|)})$.
Define $t:=
x_0^{(i)}-x^{(i)}_{\beta_i+1}$.
Then $t \cdot
(f^{>\alpha<})^{>\gamma<}(x)$
is given by
\begin{eqnarray*}
\!\!\! & &\!\!\!
(f^{>\alpha<})^{>\beta<}(
x^{(1)};\ldots; x^{(i-1)};
x^{(i)}_0, x^{(i)}_1,\ldots, x^{(i)}_{\beta_i};
x^{(i+1)};\ldots; x^{(d+|\alpha|)})\\
\!\!\! &&\!\!\! -\,
(f^{>\alpha<})^{>\beta<}(x^{(1)};\ldots; x^{(i-1)};
x^{(i)}_{\beta_i+1},x^{(i)}_1,\ldots, x^{(i)}_{\beta_i};
x^{(i+1)};\ldots; x^{(d+|\alpha|)})\\
\!\!\!&=&
\!\!\! f^{>\alpha+\bar{\beta}<}(x^{(1)};\ldots; x^{(i-1)};
x^{(i)}_0,x^{(i)}_1,\ldots, x^{(i)}_{\beta_i};
x^{(i+1)};\ldots; x^{(d+|\alpha|)})\\
\!\!\! &&\!\!\! -\,
f^{>\alpha+\bar{\beta}<} (x^{(1)};\ldots; x^{(i-1)};
x^{(i)}_{\beta_i+1},x^{(i)}_1,\ldots, x^{(i)}_{\beta_i};
x^{(i+1)};\ldots; x^{(d+|\alpha|)})\\
\!\!\!&=& \!\!\!
f^{>\alpha+\bar{\beta}<}(y; x^{(s_j)};\ldots; x^{(i-1)};
x^{(i)}_0,x^{(i)}_1,\ldots, x^{(i)}_{\beta_i};
x^{(i+1)},\ldots, x^{(s_{j+1}-1)}; z)\\
\!\!\! &&\!\!\! -\,
f^{>\alpha+\bar{\beta}<}(y; x^{(s_j)};\ldots; x^{(i-1)};
x^{(i)}_{\beta_i+1},x^{(i)}_1,\ldots, x^{(i)}_{\beta_i};
x^{(i+1)};\ldots;  x^{(s_{j+1}-1)}; z)\\
\!\!\!&=&\!\!\!
f^{>\alpha+\bar{\beta}<}(y;
x^{(i)}_0,x^{(i)}_1,\ldots, x^{(i)}_{\beta_i};
x^{(s_j)};\ldots;x^{(i-1)};
x^{(i+1)};\ldots;x^{(s_{j+1}-1)}; z)\\
\!\!\! &&\!\!\! -\,
f^{>\alpha+\bar{\beta}<} (y;
x^{(i)}_{\beta_i+1},x^{(i)}_1,\ldots, x^{(i)}_{\beta_i};
x^{(s_j)};\ldots;x^{(i-1)};
x^{(i+1)};\ldots;x^{(s_{j+1}-1)}; z)\\
\!\!\!&=&\!\!\!\!
t f^{>\alpha+\bar{\gamma}<}(y;
x^{(i)}_0,x^{(i)}_1,\ldots, x^{(i)}_{\beta_i};
x^{(s_j)};\ldots;x^{(i-1)};
x^{(i+1)};\ldots;x^{(s_{j+1}-1)}; x^{(i)}_{\beta_i+1}; z)\\
\!\!\!&=&\!\!\!
t f^{>\alpha+\bar{\gamma}<}(y;x^{(s_j)};\ldots;
x^{(s_{j+1}-1)};z)\;=\;
t f^{>\alpha+\bar{\gamma}<}(x)\,,
\end{eqnarray*}
using~(\ref{simpfml2}) for the first equality
and the symmetry properties
of $f^{>\alpha+\bar{\beta}<}$
and $f^{>\alpha+\bar{\gamma}<}$
(as in Lemma~\ref{presymm})
for the third and penultimate
equality.
This completes the inductive proof of (\ref{simpfml2}).\\[2.5mm]
As a consequence of (\ref{simpfml2}),
the continuous (resp., $C^{0,\sigma}$-) map
$f^{<\alpha+\bar{\beta}>}$
extends the map $(f^{>\alpha<})^{>\beta<}$.
It hence also extends $(f^{<\alpha>})^{>\beta<}$,
for each $\beta\in \N_0^{d+|\alpha|}$
with $|\beta|\leq k-|\alpha|$.
Hence $f^{<\alpha>}$ is $C^{k-|\alpha|}_{SDS}$
(resp., $C^{k-|\alpha|,\sigma}_{SDS}$)
and (\ref{simpfml}) holds.\vspace{2.5mm}\Punkt

\noindent
{\bf Proof of Theorem~\ref{mainthm},
completed.}
It remains to show that if $f$ is $C^k_{SDS}$
(resp., $C^{k,\sigma}_{SDS}$ with $\sigma\in \;]0,1]$),
then $f$ is $C^k_{BGN}$ (resp., $C^{k,\sigma}_{BGN}$).
We assume first that
$U=U_1\times\cdots\times U_d$
for certain subsets $U_i\sub \K$.
The proof is by induction
on $k\in \N_0$. The case $k=0$ being trivial,
assume now that $f$ is $C^k_{SDS}$ (resp.,
$C^{\,k,\sigma}_{SDS}$) for some $k\geq 1$.
For each $(x,y,t)\in U^{]1[}$,
we have
%
\begin{eqnarray}
\hspace*{-1mm}\lefteqn{f^{]1[}(x,y,t)}\quad\notag\\
& \!\! = \!\! & \frac{f(x+ty)-f(x)}{t}
\, =\,  \sum_{j=1}^d
\frac{f(x+t\sum_{i=1}^jy_ie_i)-f(x+t\sum_{i=1}^{j-1}y_ie_i)}{t}\notag\\
& \!\! = \!\! &  \sum_{j=1}^d
y_jf^{<e_j>}(x_1\! +\! ty_1;\ldots;x_{j-1}\! +\! ty_{j-1};
x_j,x_j\! + \! ty_j;x_{j+1},\ldots, x_d),\;  \label{rhs}
\end{eqnarray}
because
$a_j:=\frac{f(x+t\sum_{i=1}^jy_ie_i)-f(x+t\sum_{i=1}^{j-1}y_ie_i)}{t}
=y_j\frac{f(x+t\sum_{i=1}^jy_ie_i)-f(x+t\sum_{i=1}^{j-1}y_ie_i)}{y_jt}$
coincides with
$b_j\! :=\! y_jf^{<e_j>}(x_1+ty_1;\ldots;x_{j-1}+ty_{j-1};
x_j,x_j+ty_j;x_{j+1},\ldots, x_d)$
if $y_j\not=0$,
while both~$a_j$ and~$b_j$
vanish if $y_j=0$.
Since the right hand side of~(\ref{rhs})
defines a continuous (resp., $C^{0,\sigma}$-) map on all of $U^{[1]}$,
we see that $f$ is $C^1_{BGN}$ (resp., $C^{\,1,\sigma}_{BGN}$) with
%
%
\begin{eqnarray}
\lefteqn{f^{[1]}(x,y,t)}\quad \notag\\
\!\!\!\!\! &\! \!\! \!\! =\!\! &
\sum_{j=1}^d y_jf^{<e_j>}(x_1 + ty_1;\ldots;
x_{j-1} +  ty_{j-1}; x_j, x_j + ty_j;x_{j+1},\ldots, x_d)\;\;\; \label{fviaphi}
\end{eqnarray}
for all $(x,y,t)\in U^{[1]}$.
Here $f^{<e_j>}$
is a $C^{k-1}_{SDS}$-map (resp.,
a $C^{\,k-1,\sigma}_{SDS}$-map)
on $U^{<e_j>}=U_1\times\cdots\times U_{j-1}\times
U_j\times U_j\times U_{j+1}\times\cdots\times U_d$
by Lemma~\ref{preSDS}
and hence a $C^{k-1}_{BGN}$-map (resp.,
a $C^{\,k-1,\sigma}_{BGN}$-map),
by induction.
Formula (\ref{fviaphi})
now shows that $f^{[1]}$ is
built up from $f^{<e_j>}$
and various smooth maps,
whence~$f^{[1]}$ is $C^{k-1}_{BGN}$ (resp.,
$C^{\,k-1,\sigma}_{BGN}$), by the Chain Rule
(resp., Lemma~\ref{chainLip}).
Hence~$f$ is $C^k_{BGN}$ (resp.,
$C^{\,k,\sigma}_{BGN}$). This finishes
the proof if $U=U_1\times\cdots\times U_d$.\\[2.5mm]
If $U$ is open but not necessarily
of the form $U_1\times\cdots\times U_d$,
then every point $x\in U$
has an open neighbourhood~$V$
of the form $V=V_1\times\cdots\times V_d$
for certain open subsets $V_1,\ldots, V_d\sub \K$.
It is clear that $f|_V$ is $C^k_{SDS}$ (resp.,
$C^{\,k,\sigma}_{SDS}$) if so is~$f$.
Thus $f|_V$ is
$C^k_{BGN}$ (resp.,
$C^{\,k,\sigma}_{BGN}$)
by the case
already settled,
and thus $f$ is
$C^k_{BGN}$ (resp.,
$C^{\,k,\sigma}_{BGN}$)
as these properties can be checked locally
(see Lemma~\ref{Crlocal}, resp., Lemma~\ref{hoelloc}).\,\vspace{2.5mm}\Punkt

\noindent
In the following, we shall often simply
refer to $C^k_{BGN}$-maps as~$C^k$-maps,
and to $C^{k,\sigma}_{BGN}$-maps as~$C^{k,\sigma}$-maps.
\section{Comparison with Ludkovsky's concepts}\label{seclud}
%
%
%
In this section, we give a definition of $C^k$-maps
following an idea of Ludkovsky,
and show that such maps need not be
$C^k_{BGN}$ (nor $C^k_{SDS}$)
in the case of ground fields
of positive characteristic.
\begin{numba}
To define $C^k$-maps in Ludkovsky's
sense, we find it useful to introduce the
following notations for $U$ an open subset
of a topological vector space~$E$
over a topological field~$\K$:
We define $\Phi_1(U):=U^{]1[}$ and
$\wb{\Phi}_1(U):=U^{[1]}$. Given an integer
$k\geq 2$, we let
$\wb{\Phi}_k(U)$ be the set of all
$(x,\xi_1,\ldots, \xi_k,t_1,\ldots, t_k)\in
U\times E^k\times \K^k$
such that
$(x,\xi_1,\ldots, \xi_{k-1},t_1,\ldots, t_{k-1})\in \wb{\Phi}_{k-1}(U)$
holds as well as
$(x+t_k\xi_k,\xi_1,\ldots, \xi_{k-1},t_1,\ldots, t_{k-1})\in
\wb{\Phi}_{k-1}(U)$.
Finally, we let $\Phi_k(U)$
be the set of all
$(x,\xi_1,\ldots, \xi_k,t_1,\ldots, t_k)\in \wb{\Phi}_k(U)$
such that $t_k\not=0$.
\end{numba}
%
%
\begin{defn}\label{deflud}
Let $E$ and $F$ be topological
vector spaces over a topological
field $\K$, and $f \colon U\to F$ be a map on an open
subset $U\sub E$.
We say that $f$ is $C^1_{Lud}$
if $f$ is $C^1_{BGN}$,
i.e., if the directional difference quotient map
$\Phi_1 (f):=f^{]1[}$
admits a continuous extension
$\wb{\Phi}_1 (f):=f^{[1]}$ to
$\wb{\Phi}_1(U)=U^{[1]}$.
Recursively, having declared
when $f$ is a $C^{k-1}_{Lud}$-map
and defined a map $\wb{\Phi}_{k-1}(f)\colon \wb{\Phi}_{k-1}(U)\to F$
in this case,
we say that $f$ is $C^k_{Lud}$
if~$f$ is $C^{k-1}_{Lud}$
and if the map $\Phi_k(f)\colon \Phi_k(U)\to F$
taking $(x,\xi_1,\ldots, \xi_k,t_1,\ldots, t_k)$ to
\[
\frac{
\wb{\Phi}_{k-1}(x+t_k\xi_k,\xi_1,\ldots, \xi_{k-1},t_1,\ldots, t_{k-1})
-
\wb{\Phi}_{k-1}(x,\xi_1,\ldots, \xi_{k-1},t_1,\ldots, t_{k-1})}{t_k}
\]
admits a continuous extension
$\wb{\Phi}_k(f)\colon \wb{\Phi}_k(U)\to F$.
We say that $f$ is $C^\infty_{Lud}$
if $f$ is $C^k_{Lud}$ for each $k\in \N$.
\end{defn}
\begin{rem}
Ludkovsky defined $C^k$-maps only
for certain ultrametric fields
of characteristic~$0$ and $E$, $F$
locally convex, but of course
the preceding definition is meaningful
in the stated generality.
Furthermore, he only required the existence
of $\wb{\Phi}_1(f)$ locally
on a neighbourhood of $(x,0,0)$ for
each given point $x\in U$
(and similarly for $\wb{\Phi}_k(f)$).
We find it more convenient to define
$\wb{\Phi}_k(f)$
globally on $\wb{\Phi}_k(U)$
(which is equivalent to the local existence).
\end{rem}
\begin{rem}
If $f$ is $C^k_{Lud}$,
we define
\[
d^{(j)}f(x,\xi_1,\ldots,\xi_j)\; :=\;
\wb{\Phi}_j(f)(x,\xi_1,\ldots,\xi_j,0,\ldots, 0)
\]
for $x\in U$, $j\in \N$ with $j\leq k$,
and $\xi_1,\ldots, \xi_j\in E$.
We also write $df(x,\xi):=d^{(1)}f(x,\xi)$.
Then $d^{(j)}f\colon U\times E^j\to F$
is continuous, being a partial map of $\wb{\Phi}_j(f)$
(i.e., a map obtained from $\wb{\Phi}_j(f)$
by fixing some of its arguments).
Furthermore,
$f^{(j)}(x):= d^{(j)}f(x,\sbull)\colon E^j\to F$
is a symmetric $j$-linear map,
by a reasoning similar to that
used to prove \cite[Lemma~4.8]{BGN}.\footnote{In
the example discussed in Theorem~\ref{countludk} below,
the $j$-linearity will be obvious.}
We remark that Ludkovsky made the $j$-linearity
of the maps $f^{(j)}(x)$ part of his
definition of a $C^k$-map;
by the preceding, this requirement is
redundant and can be omitted. 
\end{rem}
\begin{rem}
We mention that Definition~\ref{deflud}
captures the basic idea of\linebreak
Ludkovsky's
approach, but differs slightly from his
actual definition which imposes
additional boundedness conditions.
In the example discussed in Theorem~\ref{countludk} below,
the domain~$\bO$ will be an open and compact
set, whence these additional
conditions will be satisfied
automatically.
\end{rem}
\begin{rem}
It is clear that each $C^k_{BGN}$-map
is also $C^k_{Lud}$; a suitable
partial map of $f^{[j]}$ serves
as the continuous extension
$\wb{\Phi}_j(f)$ of $\Phi_j(f)$,
for each $j\in \N$ such that $j\leq k$.
To make this more precise,
let us write
$E^{[j]}=E\times H_j\times \K$,
where $H_j$ collects all factors
in the middle.
Explicitly, we have
$H_1:=E$,
$H_j:= H_{j-1}\times \K\times E^{[j-1]}$
if $j\geq 2$.
Let $0_{H_j}$ be the zero element
in~$H_j$.
Then a simple induction on $k\in \N$ shows
that, if $f$ is $C^k_{BGN}$,
then $f$ is $C^k_{Lud}$, with
\[
\wb{\Phi}_k(x,\xi_1,\ldots,\xi_k,t_1,\ldots,t_k)\, :=\,
f^{[k]}(x,\xi_1,t_1;\xi_2,0_{H_1},t_2;\ldots;
\xi_k,0_{H_{k-1}},t_k)
\]
for $(x,\xi_1,\ldots,\xi_k,t_1,\ldots, \xi_k)\in \wb{\Phi}_k(U)$
providing the continuous extension
of $\Phi_k(f)$ to a map on $\wb{\Phi}_k(U)$.
\end{rem}
Let $\K$ be a local field of positive
characteristic now.
Thus, up to isomorphism,
$\K=\F_q(\!(X)\!)$ is a field of formal
Laurent series over a finite field $\F_q$
with $q=p^\ell$ elements
for some $\ell\in \N$.
We let $\bO:=\F_q[\![X]\!]$
be the ring of formal power series,
which is an open, compact subring
of~$\K$. Its elements
are of the form $x=\sum_{k=0}^\infty a_k X^k$,
where $a_k\in \F_q$.
Recall that if $x\not=0$, then
its absolute value is given by $|x|=q^{-k}$,
where $k\in \N_0$ is chosen minimal
such that $a_k\not=0$.
We define a mapping
%
\begin{equation}\label{defbd}
f\colon \bO\to \K\, ,
\quad \sum_{k=0}^\infty a_kX^k
\mto \sum_{k=0}^\infty a_kX^{[\frac{3}{2}k]}\, ,
\end{equation}
where $[r]$ denotes the Gau\ss{} bracket
(integer part) of a real number $r\geq 0$.
%
\begin{thm}\label{countludk}
The map $f\colon \bO\to \K$
defined in {\rm(\ref{defbd})}
is $C^\infty_{Lud}$,
but not $C^2_{BGN}$.
\end{thm}
\begin{proof}
It is useful to note first that $f$ is a homomorphism
of additive groups, i.e., $f(x+y)=f(x)+f(y)$
for all $x,y\in \bO$.
For $x,y\in \bO$, we have
\[
|x-y|^{\frac{3}{2}}\; \leq \;
|f(x-y)|\;\leq\; 
q |x-y|^{\frac{3}{2}}\, ,
\]
where $|f(x-y)|=|f(x)-f(y)|$.
As a consequence, $f$ is $C^1_{BGN}$
with derivative $f'(x)=0$ for all $x\in \bO$
(see \cite[Lemma 2.1]{NOA};
cf.\ \cite[Theorem 29.12]{Sch}).
Furthermore, $f$ is not $C^2_{BGN}$
because it does not admit a second order
Taylor expansion (see \cite[Lemma 2.2]{NOA}).
We now show that $f$ is $C^\infty_{Lud}$.
First, we note that $f$ is $C^1_{Lud}$
because it is $C^1_{BGN}$,
with
%
\begin{equation}\label{simsh}
\wb{\Phi}_1(f)(x,\xi,0)\, =\, df(x,\xi)\, =\, f'(x)\xi \, =\, 0
\quad \mbox{for all $x\in \bO$ and $\xi\in \K$.}
\end{equation}
Using that $f$ is a homomorphism,
for all $x\in \bO$, $\xi\in \K$
and $t\in \K^\times$ such that $x+t\xi \in \bO$,
we obtain
%
\begin{equation}\label{essl}
\wb{\Phi}_1(f)(x,\xi,t)
= \frac{f(x+t\xi)-f(x)}{t}
= \frac{f(x)+f(t\xi)-f(x)}{t}
= \frac{f(t\xi)}{t} .
\end{equation}
Since the right hand side of (\ref{essl})
is independent of~$x$, we obtain
for all $x\in \bO$, $\xi_1,\xi_2\in \K$
and $t_1,t_2\in \K^\times$
such that $(x,\xi_1,\xi_2,t_1,t_2)\in \Phi_2(\bO)$:
%
\begin{equation}\label{essl2}
\Phi_2(f)(x,\xi_1,\xi_2,t_1,t_2)=
\frac{\wb{\Phi}_1(f)(x+t_2\xi_2,\xi_1,t_1)-\wb{\Phi}_1(f)(x,\xi_1,t_1)}{t_2}
= 0 \, .
\end{equation}
By (\ref{simsh}), we also have
$\Phi_2(f)(x,\xi_1,\xi_2,t_1,t_2)=0$
for all $x\in \bO$, $\xi_1,\xi_2\in \K$,
$t_1=0$ and $t_2\in \K^\times$
such that $(x,\xi_1,\xi_2,0,t_2)\in \Phi_2(\bO)$.
Thus
\[
\wb{\Phi}_2(f)\colon \wb{\Phi}_2(\bO)\to F\, ,
\quad (x,\xi_1,\xi_2,t_1,t_2)\mto 0
\]
is a continuous map which extends
$\Phi_2(f)$, and thus $f$ is
$C^2_{Lud}$ with $\wb{\Phi}_2(f)=0$.
It now readily follows by induction
that $f$ is $C^k_{Lud}$ for each
$k\geq 2$, with $\wb{\Phi}_k(f)=0$.
\end{proof}
\begin{rem}
It would be interesting to clarify whether
$C^k_{Lud}$-maps between locally convex
spaces over an ultrametric field~$\K$ of characteristic~$0$
(as originally defined by Ludkovsky)
coincide with $C^k_{BGN}$-maps (for $k\geq 2$),
notably for $\K=\Q_p$.
It is also unknown whether
$C^k_{Lud}$-maps into real
non-locally convex spaces
are $C^k_{BGN}$.
The author conjectures that neither
is the case, but has not found counterexamples
so far.
\end{rem}
\section{H\"{o}lder differentiable maps on
metrizable spaces}\label{secmetr}
%
%
%
%
This section is devoted
to the proof of the following
characterization of
$C^{\ell,\sigma}_{BGN}$-maps
on open subsets of metrizable spaces.
%
%
\begin{thm}\label{dfc}
Let $(\K,|.|)$ be $\R$ or
an ultrametric field.
Let $E$ and $F$ be topological $\K$-vector spaces
and $f\colon U\to F$ be a map,
defined on an open subset $U\sub E$.
Let $\ell\in \N_0$ and $\sigma\in \; ]0,1]$.
If $E$ is metrizable,
then $f$ is $C^{\,\ell,\sigma}$
if and only if $f\circ \gamma\colon \K^{\ell+1}\to F$
is $C^{\,\ell,\sigma}$, for each
$C^\infty$-map
$\gamma \colon \K^{\ell+1}\to U$.
\end{thm}
The proof of Theorem~\ref{dfc}
heavily relies on tools
developed in~\cite{GaLa},
which are variants of standard methods
of differential calculus in
real locally convex spaces (cf.\ \cite{KaM}).
To describe these tools,
we need the auxiliary notion
of a ``calibration''
on a topological vector space~$E$
over a valued field~$\K$.
%
%
\begin{defn}\label{defcalib}
A sequence $(q_n)_{n\in \N_0}$
of gauges on~$E$ is called a \emph{calibration} if
%
\begin{equation}\label{ptycalib}
(\forall n\in\N_0)(\forall x,y\in E)\quad
q_n(x+y)\;\leq\; q_{n+1}(x)+q_{n+1}(y)\,.
\end{equation}
If $q$ is a gauge on~$E$, then there always
exists a calibration $(q_n)_{n\in \N_0}$
such that $q_0=q$ (cf.\ Remark~\ref{remfk});
we then say that $q$ extends to $(q_n)_{n\in \N_0}$.
\end{defn}
%
%
\begin{rem}\label{increas}\label{excalib}
If $(q_n)_{n\in \N_0}$ is a calibration,
then $q_n\leq q_{n+1}$ for each $n\in \N_0$
because $q_n(x)=q_n(x+0)\leq q_{n+1}(x)+q_{n+1}(0)=q_{n+1}(x)$
for each $x\in E$.
Also note that if $q\colon E\to [0,\infty[$ is a continuous
seminorm, then $(q)_{n\in \N_0}$
is a calibration.
If $(q_n)_{n\in \N_0}$
is any calibration extending the seminorm~$q$,
then $q_n\geq q$ for each $n$,
by the preceding remark.
Thus $(q)_{n\in \N_0}$ is the smallest
calibration extending~$q$.
\end{rem}
The following two lemmas are
the main results of~\cite{GaLa}.
They are variants of \cite[Lemma~12.2]{KaM}.
In the first lemma, $\bO:=\wb{B}^\K_1(0)$.
%
%
%
\begin{la}[Ultrametric General Curve Lemma]\label{thmgencu}
Let $E$ be a topological vector space
over an ultrametric field~$\K$,
$\rho\in \K^\times$ with $|\rho|<1$
and $(\gamma_n)_{n\in \N}$ be a family
of smooth maps $\gamma_n\in BC^\infty(\rho^n\bO, E)$
which become small sufficiently fast
in the sense that,
for each gauge $q$ on~$E$,
there exists a calibration $(q_n)_{n\in\N_0}$
extending $q$ such that
%
\begin{equation}\label{condgen}
(\forall a>0)\, (\forall k, m\in \N_0)\quad
\lim_{n\to\infty}\,
a^n\|\gamma_n^{<k>}\|_{q_{n+m}, \infty}\;=\;  0\, .
\end{equation}
Then there exists a smooth curve
$\gamma\colon \K\to E$ with
$\im(\gamma)=\{0\}\cup\bigcup_{n\in \N}\im(\gamma_n)$,
such that
$\gamma(\rho^{n-1}+t)=
\gamma_n(t)$
for all $n\in \N$ and $t\in \rho^n\bO$.\,\Punkt
\end{la}
%
%
\begin{rem}\label{symetr}
Let $E$ in Lemma~\ref{thmgencu} be metrizable
and suppose that there exists
a calibration $(p_n)_{n\in \N_0}$
such that $\{p_n\colon n\in \N_0\}$
is a fundamental system of gauges,
and $C>0$ such that
%
\begin{equation}\label{shalluse}
(\forall k\in \N_0)\, (\forall n\geq k)\quad
\|\gamma_n^{<k>}\|_{p_{2n},\infty}\;\leq\; C n^{-n}\,.
\end{equation}
Then the hypothesis (\ref{condgen})
of Lemma~\ref{thmgencu} is satisfied:
Given~$q$, we can extend it to a suitable calibration
via $q_n:=rp_{n+n_0}$ for $n \in \N$,
with $r>0$ and $n_0\in \N_0$ sufficiently large.
\end{rem}
%
%
\begin{la}[Real Case of General Curve Lemma]\label{thmrealcu}
Let $E$ be a real to\-pological
vector space and $(s_n)_{n\in \N}$
as well as $(r_n)_{n\in \N}$
be sequences of positive reals such that
$\sum_{n=1}^\infty s_n<\infty$ and $r_n\geq
s_n+\frac{2}{n^2}$ for each $n\in \N$. Furthermore,
let $(\gamma_n)_{n\in \N}$ be a sequence of smooth
maps $\gamma_n\colon [{-r_n}, r_n] \to E$
which become small sufficiently
fast in the sense that, for each gauge~$q$ on~$E$,
there exists a calibration $(q_n)_{n\in \N_0}$
extending $q$ such that
%
\begin{equation}\label{dfco}
(\forall k,\ell, m\in \N_0)\qquad
\lim_{n\to\infty} \, n^\ell\|\gamma_n^{<k>}\|_{q_{n+m},\infty}\,=\,0\,.
\end{equation}
Then there exists a curve $\gamma\in BC^\infty(\R, E)$
with $\,\im(\gamma)\sub [0,1]\cdot\bigcup_{n\in \N}\im(\gamma_n)$
and a convergent sequence $(t_n)_{n\in \N}$
of real numbers such that
$\gamma(t_n+t) = \gamma_n(t)$
for all $ n\in \N$ and $ t\in [-s_n,s_n]$.\Punkt
\end{la}
Again, (\ref{shalluse}) enables
to manufacture calibrations
satisfying~(\ref{dfco}).\\[2.5mm]
In our applications, the maps $\gamma_n$
are restrictions
of affine-linear maps to balls.
The following simple lemma will help us to
verify the hypotheses of the General
Curve Lemmas in this case.
%
%
%
%
\begin{la}\label{deg1}
Consider the map $\gamma\colon \wb{B}_r^\K(0)\to E$,
$x\mto xa+b$,
where
$E$ is a topological vector space
over a valued field~$\K$,
$a,b\in E$ and $r>0$.
Let $q_1$ and $q_2$ be gauges
on~$E$ such that $q_1(x+y)\leq q_2(x)+q_2(y)$
for all $x,y\in E$.
Then
$\|\gamma\|_{q_1,\infty}\leq r\|a\|_{q_2}+\|b\|_{q_2}$,
$\|\gamma^{<1>}\|_{q_1,\infty}=\|a\|_{q_1}$
and $\|\gamma^{<k>}\|_{q_1,\infty}=0$
for $k\geq 2$.
\end{la}
\begin{proof}
Since $\|\gamma(x)\|_{q_1}=\|xa+b\|_{q_1}\leq
|x|\cdot\|a\|_{q_2}+\|b\|_{q_2}\leq
r\cdot\|a\|_{q_2}+\|b\|_{q_2}$ for each $x\in \wb{B}_r^\K(0)$,
the first inequality holds.
The remaining assertions follow
from the observations that $\gamma^{<1>}(x,y)=a$
for all $x,y\in \wb{B}_r^\K(0)$
and $\gamma^{<k>}=0$
for all $k\geq 2$.
\end{proof}
Another simple observation will be used.
%
\begin{la}\label{ctoff}
Let $(\K,|.|)$ be
either $\R$ or
an ultrametric field.
Let $E$ and $F$ be topological $\K$-vector spaces
and $f\colon U\to F$ be a map,
defined on an open subset $U\sub E$.
Let $\ell, d \in \N_0$ and $\sigma\in \;]0,1]$.
If $f\circ \gamma\colon \K^d\to F$
is $C^{\,\ell,\sigma}$, for each
$C^\infty$-map
$\gamma \colon \K^d \to U$,
then also
$f\circ \gamma\colon  V \to F$
is $C^{\,\ell,\sigma}$, for each
$C^\infty$-map
$\gamma \colon V\to U$
defined on an open subset $V\sub \K^d$.
\end{la}
\begin{proof}
Given $x_0\in V$, there exists a smooth map
$\kappa\colon \K^d\to V$
such that $\kappa|_W=\id_W$
for some open neighbourhood $W\sub V$
of~$x_0$.
In fact, if~$\K$ is ultrametric,
we can choose an open, closed neighbourhood
$W\sub V$ of~$x_0$ and define $\kappa(x):=x$
if $x\in W$, $\kappa(x):=x_0$ if $x\in V\setminus W$.
In the real case, we can manufacture~$\kappa$
by standard arguments, using a cut-off
function. Then $\eta:=\gamma\circ \kappa\colon \K^d\to U$
is smooth and hence
$f\circ \eta$
is $C^{\ell,\sigma}$.
Then $(f\circ\gamma)|_W=
(f\circ \eta)|_W$ is~$C^{\ell,\sigma}$.
Hence $f\circ\gamma$ is locally $C^{\ell,\sigma}$
and thus $C^{\ell,\sigma}$, by Lemma~\ref{hoelloc}.
\end{proof}
{\bf Proof of Theorem~\ref{dfc}.}
If $f$ is $C^{\ell,\sigma}$, then $f\circ \gamma$ is $C^{\ell,\sigma}$
for each $C^\infty$-map $\gamma \colon \K^{\ell+1} \to U$,
by Lemma~\ref{chainLip}
and Remark~\ref{autolip}.
To prove the converse direction,
we first assume that $\K$ is an ultrametric field.
We start with the case $\ell=0$.
If~$f$ is not $C^{0,\sigma}$,
then the condition formulated
in Definition~\ref{deflipcts}
is violated by some $x_0\in U$.
Hence, there exists a gauge~$q$ on~$F$ such that,
for each neighbourhood $V\sub U$
of~$x_0$ and gauge~$p$ on~$E$,
there are $x,y\in V$ such that
$\|f(x)-f(y)\|_q>(\|x-y\|_p)^\sigma$.
After a translation, we may assume
that $x_0=0$.
Pick a gauge $q_0$ on~$E$
such that $B_1^{q_0}(0)\sub U$,
and extend it to a calibration
$(q_n)_{n\in \N_0}$
on~$E$
such that $\{q_n\colon n\in \N_0\}$
is a fundamental system of gauges.
Also, pick $\rho\in \K^\times$
such that $|\rho|<1$.
After replacing $q_1, q_2,\ldots$
by large multiples if necessary, we
may assume that
%
\begin{equation}\label{preprp}
q_0 \;\leq\;
\big({\textstyle
\frac{1}{2}+\frac{1}{|\rho|}}
\big)^{-1}  q_1\, .
\end{equation}
Applying the above property of~$q$ for a given $n\in \N$
to $V:=\wb{B}^{q_{2n+3}}_{\frac{1}{2}n^{-n}|\rho|^n}(0)$
and $p:=n^{\frac{1}{\sigma}}n^n q_{2n+2}$,
we find $x_n,y_n\in E$ such that
%
%
\begin{equation}\label{firstprp}
\|x_n\|_{q_{2n+3}},\|y_n\|_{q_{2n+3}}
\;\leq\; \frac{1}{2}n^{-n}|\rho|^n
\end{equation}
and
%
\begin{equation}\label{secprp}
\|f(x_n)-f(y_n)\|_q\; >\;
n\cdot n^{\sigma n}\big(\|x_n-y_n\|_{q_{2n+2}}\big)^\sigma\,.
\end{equation}
\emph{Case}~1. If $\|x_n-y_n\|_{q_{2n+2}}\not=0$,
let $k_n$ be the unique integer such that
%
\begin{equation}\label{stasta}
|\rho|^{k_n}\; \leq \; n^n\|x_n-y_n\|_{q_{2n+2}}\; < \; |\rho|^{k_n-1}\, .
\end{equation}
Since $n^n\|x_n-y_n\|_{q_{2n+2}}
\leq n^n(\|x_n\|_{q_{2n+3}}+\|y_n\|_{q_{2n+3}})\leq
|\rho|^n$ by (\ref{firstprp}),
we have $k_n\geq n$.\\[3mm]
\emph{Case}~2. If $\|x_n-y_n\|_{q_{2n+2}}=0$,
we choose the integer $k_n\geq n$ so large that
$\|f(x_n)-f(y_n)\|_q\geq n (|\rho|^{k_n})^\sigma$.\\[3mm]
In either case, we define
\[
\gamma_n\colon \wb{B}^{\K}_{|\rho|^n}(0)\to E\,,\quad
\gamma_n(t)\, :=\, x_n+\frac{t}{\rho^{k_n}}(y_n-x_n)\,.
\]
By Lemma~\ref{deg1}, we then have
$\gamma_n^{<k>}=0$ for $k\geq 2$,
furthermore
\[
\|\gamma_n^{<1>}\|_{q_{2n},\infty}
\;=\;
\frac{\|x_n-y_n\|_{q_{2n}}}{|\rho|^{k_n}}
\;\leq\;
\frac{\|x_n-y_n\|_{q_{2n+2}}}{|\rho|^{k_n}}
\;<\; \frac{n^{-n}}{|\rho|}
\]
by definition of $k_n$,
and finally
\[
\|\gamma_n\|_{q_{2n},\infty}\leq
|\rho^n|\frac{\|x_n-y_n\|_{q_{2n+1}}}{|\rho|^{k_n}}
+\|x_n\|_{q_{2n+1}}
<
{\textstyle \frac{n^{-n}}{|\rho|}+\frac{1}{2}n^{-n}}|\rho|^n
< {\textstyle \big(\frac{1}{2}+\frac{1}{|\rho|}\big)}n^{-n}\, ,
\]
entailing that $\|\gamma_n \|_{q_0,\infty}<1$
(see (\ref{preprp})) and thus
$\,\im \gamma_n\sub B^{q_0}_1(0)\sub U$.
In view of the preceding,
(\ref{shalluse}) in Remark~\ref{symetr}
is satisfied by the calibration $(q_n)_{n\in \N_0}$
with $C=\frac{1}{2}+\frac{1}{|\rho|}$.
Therefore the General Curve Lemma (Lemma~\ref{thmgencu})
provides a smooth map $\gamma\colon \K\to E$
with $\gamma(\K)\sub U$ such that
$\gamma(t)=\gamma_n(t-\rho^{n-1})$
for each $n\in \N$ and $t\in \wb{B}_{|\rho|^n}^\K(\rho^{n-1})$.
In particular,
$\gamma(\rho^{n-1})=\gamma_n(0)=x_n$
and $\gamma(\rho^{n-1}+\rho^{k_n})=\gamma_n(\rho^{k_n})
=y_n$ for each $n\in \N$.
Hence
%
%
\begin{eqnarray}
\hspace*{-2cm}
\lefteqn{\|f(\gamma(\rho^n))-f(\gamma(\rho^{n-1} +\rho^{k_n})) \|_q}
\qquad \notag \\[.5mm]
&=& \|f(x_n)-f(y_n)\|_q\;>\;
n\cdot n^{\sigma n} (\|x_n-y_n\|_{q_{2n+2}})^\sigma \label{sta}\\
& = & 
n\cdot n^{\sigma k}\left(
\frac{\|x_n-y_n\|_{q_{2n+2}}}{|\rho^{k_n}|} \right)^\sigma
|\rho^{k_n}|^\sigma
\;\geq\; n |\rho^{k_n}|^\sigma\label{oneref}
\end{eqnarray}
in Case~1, using (\ref{secprp})
in (\ref{sta}) and then (\ref{stasta})
for the final inequality.
In Case~2, we have
%
\begin{equation}\label{stp}
\|f(\gamma(\rho^n))-f(\gamma(\rho^{n-1}+\rho^{k_n}))\|_q
\,=\, \|f(x_n)-f(y_n)\|_q\, \geq \, n\cdot (|\rho|^{k_n})^\sigma
\end{equation}
by choice of~$k_n$.
Thus (\ref{stp})
holds for each~$n$,
and hence $f\circ\gamma$ is not $C^{0,\sigma}$.
In fact, if $f\circ  \gamma$ were $C^{0,\sigma}$,
there would be a $0$-neighbourhood $J\sub \K$
and a gauge $g$ on $\K$
such that
$\|f(\gamma(t))-f(\gamma(s))\|_q\leq
(\|t-s\|_g)^\sigma$ for all $s,t\in J$.
As a consequence of Lemma \ref{findfund},
there is $C>0$ such that
$g\leq C|.|$. Hence
$\|f(\gamma(t))-f(\gamma(s))\|_q\leq
C^\sigma |t-s|^\sigma$ for all $s,t\in J$,
which contradicts (\ref{stp}).\\[3mm]
The general case: Let $\ell$
be a positive integer now.
If $f\circ\gamma$ is $C^{\ell,\sigma}$
for each smooth map
$\gamma\colon \K^{\ell+1}\to U$,
then
$f\circ\gamma$ is $C^\ell$
in particular and hence $f$ is $C^\ell$,
by \cite[Theorem~12.4]{BGN}.
To prove that $f$ is
$C^{\ell,\sigma}$,
it only remains to show
that $f^{[\ell]}$ is
$C^{0,\sigma}$.
We assume
that $f^{[\ell]}$ is
not $C^{0,\sigma}$
and derive a contradiction.
Since $U^{]\ell[}$
is dense in the domain
$U^{[\ell]}$ of the continuous map
$f^{[\ell]}$,
Lemma~\ref{densehoel} shows that
there exists $x_0\in U^{[\ell]}$ and
a gauge~$q$ on~$F$ such that,
for each neighbourhood $V\sub U^{[\ell]}$
of~$x_0$ and gauge~$p$ on~$E^{[\ell]}$,
there are $x,y\in V\cap U^{]\ell[}$ such that
$\|f^{[\ell]}(x)-f^{[\ell]}(y)\|_q>(\|x-y\|_p)^\sigma$.
We now pick
$x_n,y_n \in U^{]\ell[}$
as above in the case $\ell=0$,
applied to $f^{[\ell]}$
instead of~$f$,
and obtain a smooth curve $\gamma\colon
\K\to U^{[\ell]}$ such that $\gamma(\rho^{n-1})=x_n$
and $\gamma(\rho^{n-1}+\rho^n)=y_n$.
Applying \cite[Lemma~12.3]{BGN}
with $m:=1$, $V:=\K$,
$D:=\{\rho^{n-1}\colon n\in \N\}
\cup \{\rho^{n-1}+\rho^n \colon n\in \N\}$
and $X_0:=\{0\}$,
we obtain a smooth map
$\Gamma\colon W\to U$,
defined on an open subset
$W\sub \K^{\ell+1}$,
an open neighbourhood $Y$ of~$0$
in~$\K$,
and a smooth map
$g\colon Y\to W^{[\ell]}$
such that
%
\begin{equation}\label{givctr}
(\forall t\in D\cap Y)\quad
f^{[\ell]}(\gamma(t))\;=\;
(f\circ\Gamma)^{[\ell]}(g(t))\,.
\end{equation}
There is $N\in \N$ such that
$\rho^{n-1}\in Y$
and $\rho^{n-1}+\rho^n \in Y$
for each integer $n\geq N$.
The hypothesis implies
that $f\circ \Gamma$
is $C^{\ell ,\sigma}$ (see Lemma~\ref{ctoff}).
As a consequence,
$(f\circ\Gamma)^{[\ell]}$
is $C^{0,\sigma}$ and hence also
$(f\circ\Gamma)^{[\ell]}\circ g$
is $C^{0,\sigma}$.
However, by construction of~$\gamma$
and (\ref{givctr}), for each $n\geq N$ we have
\begin{eqnarray*}
\lefteqn{\|(f\circ\Gamma)^{[\ell]}(g(\rho^{n-1}))-
(f\circ\Gamma)^{[\ell]}(g(\rho^{n-1}+\rho^n))\|_q}\quad\\
&=&
\|f^{[\ell]}(\gamma(\rho^{n-1}))
- f^{[\ell]}(\gamma(\rho^{n-1}+\rho^n))\|_q\; =\;
\|f^{[\ell]}(x_n)
- f^{[\ell]}(y_n)\|_q\\
& \geq & n |\rho^n|^\sigma\,,
\end{eqnarray*}
arguing as in (\ref{oneref}) to pass to the last
line.
Hence $(f\circ\Gamma)^{[\ell]}\circ g$
is not $C^{0,\sigma}$, contradicting
the preceding. This closes
the proof in the ultrametric case.\\[3mm]
Now assume that $\K=\R$,
and pick $r\in \;]0,1[$.
If~$f$ is not $C^{0,\sigma}$,
then there exists $x_0\in U$ and a gauge~$q$ on~$F$ such that,
for each neighbourhood $V\sub U$
of~$x_0$ and gauge~$p$ on~$E$,
there are $x,y\in V$ such that
$\|f(x)-f(y)\|_q>(\|x-y\|_p)^\sigma$.
After a translation, we may assume
that $x_0=0$.
Take a gauge~$q_0$ on~$E$
such that $\wb{B}_1^{q_0}(0)\sub U$,
and extend it to a calibration
$(q_n)_{n\in \N_0}$
on~$E$
such that $\{q_n\colon n\in \N_0\}$
is a fundamental system of gauges.
We may assume that
$q_0 \leq \frac{2}{7} q_1$.
Applying the above property of~$q$ for a given $n\in \N$
to $V:=\wb{B}^{q_{2n+3}}_{\frac{1}{2}n^{-n}r^n}(0)$
and $p:=n^{\frac{1}{\sigma}}n^n q_{2n+2}$,
we find $x_n,y_n\in E$ such that
\[
\|x_n\|_{q_{2n+3}},\|y_n\|_{q_{2n+3}}
\;\leq\; \frac{1}{2}n^{-n} r^n
\]
and
$\|f(x_n)-f(y_n)\|_q >
n\cdot n^{\sigma n}\big(\|x_n-y_n\|_{q_{2n+2}}\big)^\sigma$.\\[2.5mm]
\emph{Case}~1: If $\|x_n-y_n\|_{q_{2n+2}}\not=0$,
define $s_n:=
n^n\|x_n-y_n\|_{q_{2n+2}}\leq r^n$.
\emph{Case}~2:
If $\|x_n-y_n\|_{q_{2n+2}}=0$,
choose
$s_n\in \; ]0, r^n]$ such that
$\|f(x_n)-f(y_n)\|_q\geq n \cdot (s_n)^\sigma$.
In either case, we define
$r_n:= s_n+\frac{2}{n^2}$ and
\[
\gamma_n\colon [{-r_n}, r_n] \to E\,,\quad
\gamma_n(t)\, :=\, x_n+\frac{t}{s_n}(y_n-x_n)\,.
\]
By Lemma~\ref{deg1}, we then have
$\gamma_n^{<k>}=0$ for $k\geq 2$,
furthermore
$\|\gamma_n^{<1>}\|_{q_{2n},\infty}
=\frac{\|y_n-x_n\|_{q_{2n}}}{s_n}
\leq n^{-n}$
by definition of $s_n$,
and finally
$\|\gamma_n\|_{q_{2n},\infty}<
\frac{7}{2} n^{-n}$ because
\[
\|\gamma_n(x)\|_{q_{2n}} \leq
\|x_n\|_{q_{2n+1}}+r_n\frac{\|y_n-x_n\|_{q_{2n+1}}}{s_n}
\leq
{\textstyle \frac{1}{2}}n^{-n}r^n+\left(r^n+{\textstyle \frac{2}{n^2}}\right)
n^{-n}\leq
{\textstyle \frac{7}{2}n^{-n}} ,
\]
entailing that $\|\gamma_n\|_{q_0,\infty}\leq 1$
and thus
$\,\im \gamma_n\sub \wb{B}^{q_0}_1(0)\sub U$.
In view of the preceding,
(\ref{shalluse}) in Remark~\ref{symetr}
is satisfied with $C=\frac{7}{2}$.
Therefore the General Curve Lemma (Lemma~\ref{thmrealcu})
provides a smooth map $\gamma\colon \R\to E$
with $\gamma(\R)\sub [0,1]\wb{B}^{q_0}_1(0)
=\wb{B}^{q_0}_1(0)\sub U$,
and a convergent sequence $(t_n)_{n\in \N}$
of reals such that
$\gamma(t_n+t)=\gamma_n(t)$
for each $n\in \N$ and $t\in \R$ such that
$|t|\leq s_n$.
In particular,
$\gamma(t_n)=\gamma_n(0)=x_n$
and $\gamma(t_n+s_n)=\gamma_n(s_n)
=y_n$ for each $n\in \N$.
Hence
\begin{eqnarray*}
\hspace*{-2cm}
\lefteqn{\|f(\gamma(t_n))-f(\gamma(t_n+s_n)) \|_q}
\qquad \\[.5mm]
&=& \|f(x_n)-f(y_n)\|_q\;>\;
n\cdot n^{\sigma n} (\|x_n-y_n\|_{q_{2n+2}})^\sigma\\
& = & 
n\cdot n^{\sigma k}\left(
\frac{\|x_n-y_n\|_{q_{2n+2}}}{s_n} \right)^\sigma
(s_n)^\sigma
\;=\; n \cdot (s_n)^\sigma
\end{eqnarray*}
in Case~1.
In Case~2, we have
%
\begin{equation}\label{stpR}
\|f(\gamma(t_n))-f(\gamma(t_n+s_n))\|_q
\, =\, \|f(x_n)-f(y_n)\|_q\, \geq \, n \cdot (s_n)^\sigma
\end{equation}
by choice of~$s_n$.
Thus (\ref{stpR})
holds for each~$n$,
and hence $f\circ\gamma$ is not $C^{0,\sigma}$.\\[3mm]
The general case: If~$\ell$
is a positive integer and
$f$ is not $C^{\ell,\sigma}$
although $f\circ\gamma$ is $C^{\ell,\sigma}$
for each smooth map
$\gamma\colon \K^{\ell+1}\to U$,
we reach a contradiction
along the lines of the ultrametric case.
First, applying the case $\ell=0$
to~$f^{[\ell]}$ instead of~$f$,
we find a gauge
$q$ on~$F$ and $x_n, y_n \in U^{]\ell[}$,
positive reals~$s_n$ such that
$\sum_{n=1}^\infty s_n<\infty$,
a smooth curve $\gamma\colon \R\to U^{[\ell]}$
and a convergent sequence $(t_n)_{n \in \N}$
of reals such that
$\gamma(t_n)=x_n$, $\gamma(t_n+s_n)=y_n$
and $\|f^{[\ell]}(x_n)-f^{[\ell]}(y_n)\|_q\geq n (s_n)^\sigma$
for each $n\in \N$.
Let $t_\infty :=\lim_{n\to\infty} t_n$.
Applying \cite[Lemma~12.3]{BGN}
with $m:=1$, $V:=\K$,
$D:=\{t_n \colon n\in \N\}
\cup \{t_n + s_n \colon n\in \N\}$
and $X_0:=\{t_\infty\}$,
we obtain a smooth map
$\Gamma\colon W\to U$,
defined on an open subset
$W\sub \R^{\ell+1}$,
an open neighbourhood $Y$ of~$t_\infty$
in~$\R$,
and a smooth map
$g\colon Y\to W^{[\ell]}$
such that~(\ref{givctr}) holds.
There is $N\in \N$
such that $t_n\in Y$ and $t_n+s_n\in Y$
for all $n\geq N$.
The hypothesis implies
that $f\circ \Gamma$
is $C^{\ell ,\sigma}$ (Lemma~\ref{ctoff}).
As a consequence,
$(f\circ\Gamma)^{[\ell]}$
is $C^{0,\sigma}$ and hence also
$(f\circ\Gamma)^{[\ell]}\circ g$
is $C^{0,\sigma}$.
However, by construction of~$\gamma$
and (\ref{givctr}), we have
$\|(f\circ\Gamma)^{[\ell]}(g(t_n))-
(f\circ\Gamma)^{[\ell]}(g(t_n+s_n))\|_q
=
\|f^{[\ell]}(\gamma(t_n))
- f^{[\ell]}(\gamma(t_n +s_n))\|_q=
\|f^{[\ell]}(x_n)
- f^{[\ell]}(y_n)\|_q
\geq n (s_n)^\sigma$
for each $n\geq N$,
whence
$(f\circ\Gamma)^{[\ell]}\circ g$
is not $C^{0,\sigma}$, which is absurd.\Punkt
\section{Weakly H\"{o}lder differentiable maps}\label{secweak}
%
%
%
If~$\K$ is a topological
field and $E$ a topological $\K$-vector space,
we let $E'$ be the space of all continuous
linear functionals $\lambda\colon E \to\K$.
\begin{defn}
Let $E$ and $F$ be topological
vector spaces over a valued field~$\K$
and $f\colon U\to F$ be 
a map on a subset $U\sub E$.
Let $\sigma>0 $.
We say that $f$ is \emph{weakly $C^{0,\sigma}$}
if $\lambda\circ f\colon U\to\K$
is $C^{0,\sigma}$ for each $\lambda\in F'$.
If $U$ has dense interior
and $k\in \N \cup\{\infty\}$,
we say that
$f$ is \emph{weakly $C^{k,\sigma}$}
if $\lambda\circ f\colon U\to\K$
is $C^{k,\sigma}$ for each $\lambda\in F'$.
\end{defn}
%
%
\begin{rem}\label{simplwea1}
Note that
each $C^{k,\sigma}$-map
is weakly $C^{k,\sigma}$
(cf.\ Lemma~\ref{chainlin}).
\end{rem}
%
%
\begin{rem}\label{simplwea2}
Let $f\colon U\to F$
be a weakly $C^{k,\sigma}$-map
on a subset $U\sub E$
with dense interior and
$g\colon V\to U$ be a $C^{k,1}$-map
on a subset $V$ with dense
interior of a topological
$\K$-vector space~$H$ (e.g.,
a $C^{k+1}$-map).
Then $\lambda\circ (f\circ g)=(\lambda\circ f)\circ g$
is $C^{k,\sigma}$ for each $\lambda\in F'$
(by Lemma~\ref{chainLip})
and thus $f\circ g$ is weakly~$C^{k,\sigma}$.
\end{rem}
Recall that a topological vector space
over an ultrametric field
is called\linebreak
\emph{locally convex}
if its vector topology can be defined
by a family of ultrametric seminorms
(cf.\ \cite{Mon} for further
information).
%
%
\begin{la}\label{cpinnormed}
Let $(E,\|.\|)$ be a normed space over a locally compact
field~$\K$,
$F$ be a locally convex space over~$\K$
and $f\colon K\to F$ be a
map on a compact set $K\sub E$.
Let $\sigma > 0 $.
Then the following
conditions are equivalent:
\begin{itemize}
\item[\rm(a)]
$f$ is a $C^{0,\sigma}$-map.
\item[\rm(b)]
For each continuous seminorm
$q$ on~$F$,
there is $C\in [0,\infty[$
such that
%
\begin{equation}\label{ordinhoel}
\|f(y)-f(x)\|_q\;\leq\; C\, (\|y-x\|)^\sigma \quad
\mbox{for all $x,y\in K$.}
\end{equation}
\item[\rm(c)]
$f$ is weakly $C^{0,\sigma}$.
\end{itemize}
\end{la}
\begin{proof}
(a)$\impl$(b):
Let $q$ be as in~(b).
If $f$ is $C^{0,\sigma}$,
then for each $z\in K$
there exists an open neighbourhood
$U_z\sub K$ of~$z$ and a gauge
$p_z$ on~$E$ such that
\[
\|f(y)-f(x)\|_q\;\leq\; (\|y-x\|_{p_z})^\sigma\quad
\mbox{for all $x,y\in U_z$.}
\]
Let $V_z$ be an open
neighbourhood of~$z$ in~$K$ with compact closure
$\wb{V_z}\sub U_z$.
There exists a finite subset
$\Phi\sub K$
such that $K=\bigcup_{z\in \Phi}V_z$.
For each $z\in \Phi$, there exists
$r_z>0$ such that $p_z\leq r_z\|.\|$
(cf.\ Lemma~\ref{findfund}).
Let $r:=\max\{r_z\colon z\in \Phi\}$.
The sets $\wb{V_z}$ and $K\setminus U_z$
being compact and disjoint, we can define
\[
s\; :=\; \sup\,\{\|y-x\|^{-\sigma} \colon
\mbox{$z\in \Phi$, $x\in \wb{V_z}$, $y\in K\setminus U_z$}\}\;
\in \; [0,\infty[\,.
\]
Then (\ref{ordinhoel}) holds
with $C :=\max\{r^\sigma , 2s\max \|f(K)\|_q\}$.
In fact, given \mbox{$x,y\in K$,}
there exists $z\in \Phi$
such that $x\in V_z $.
If $y\in U_z$, then
$\|f(y)-f(x)\|_q\leq$\linebreak
$(\|y-x\|_{p_z})^\sigma
\leq (r_z)^\sigma\|y-x\|^\sigma
\leq C \|y-x\|^\sigma$.
If $y\not\in U_x$, then
$\|f(y)-f(x)\|_q\leq \frac{\|f(y)\|_q+\|f(x)\|_q}{\|y-x\|^\sigma}
\|y-x\|^\sigma\leq 2s\max\|f(K)\|_q \|y-x\|^\sigma\leq
C\|y-x\|^\sigma$ as well.\vspace{1.3mm}

(b)$\impl$(a): Given a gauge $g$ on~$F$,
by local convexity there exists a continuous
seminorm~$q$ (which can be chosen ultrametric if~$\K$
is a local field) such that
$g \leq q$ (cf.\ Lemma \ref{findfund}).
Let~$C$ be as
in~(b). Then $p:=C^{\frac{1}{\sigma}}\|.\|$
is a gauge on~$E$ such that
$\|f(y)-f(x)\|_g\leq \|f(y)-f(x)\|_q\leq (\|x-y\|_p)^\sigma$
for all $x,y\in K$. Thus (a) holds.\vspace{1.3mm}

(a)$\impl$(c): See Remark~\ref{simplwea1}.\vspace{1.3mm}

(c)$\impl$(b): Pick $\rho\in \K^\times$
with $|\rho|<1$ and
define $\beta(r):=\rho^k$
for $r\in \;]0,\infty[$,
where $k\in \Z$ is the unique
integer such that $|\rho|^{k+1}< r^\sigma\leq
|\rho|^k$.
Then
%
\begin{equation}\label{estbet}
|\rho|\cdot |\beta(r)|\; < \; r^\sigma\;\leq\;
|\beta(r)|\quad\mbox{for all $r\in \;]0,\infty[$.}
\end{equation}
If $f$ is weakly $C^{0,\sigma}$, define
\[
B\;:=\;
\left\{\frac{f(y)-f(x)}{\beta(\|y-x\|)} \colon
\mbox{$x,y\in K$ such that $x\not=y$.}\right\}
\]
We claim that $B$ is bounded.
If this is true, then
$M:=\sup \|B\|_q<\infty$ for each
continuous seminorm (or gauge) $q$ on~$F$
(see Lemma~\ref{bdviagau})
and hence
\begin{eqnarray*}
\|f(y)-f(x)\|_q
& = &
\frac{\|f(y)-f(x)\|_q}{|\beta(\|y-x\|)|} \, |\beta(\|y-x\|)|
\; \leq\; M\,|\beta(\|y-x\|)|\\
& \leq &M |\rho|^{-1} \|y-x\|^\sigma
\end{eqnarray*}
for all $x,y\in K$ such that $x\not=y$,
using (\ref{estbet})
for the final inequality.
Hence (\ref{ordinhoel}) holds
with $C:=M |\rho|^{-1}$.\\[3mm]
It remains to show that~$B$ is bounded,
or equivalently, that $\lambda(B)\sub \K$
is bounded for each $\lambda\in F'$
(see \cite[Theorem~3.18]{Rud} for the real
case (from which the complex case
follows) and
\cite[Theorem~4.21]{Thi}
for the case where~$\K$ is a local
field). However, for each
$\lambda\in F'$,
the map $\lambda\circ f\colon K\to\K$ is
$C^{0,\sigma}$
and hence, by (a)$\impl$(b)
already established, there exists $C\in [0,\infty[$
such that
\[
|\lambda(f(y))-\lambda(f(x))|\;\leq\;
C\,\|y-x\|^\sigma\quad\mbox{for all $x,y\in K$.}
\]
But then $\sup|\lambda(B)|\leq C$
(whence $\lambda(B)$ is bounded),
since
%
\begin{equation}\label{mayb}
\left| \lambda\Big(\frac{f(y)-f(x)}{\beta(\|y-x\|)}\Big)\right|
=
\frac{|\lambda(f(y))-\lambda(f(x))|}{|\beta(\|y-x\|)|}
\leq
\frac{|\lambda(f(y))-\lambda(f(x))|}{\|y-x\|^\sigma} \leq
C
\end{equation}
for all $x,y\in K$ such that $x\not=y$,
using (\ref{estbet}) to obtain the first inequality.
\end{proof}
%
%
\begin{rem}\label{onlyult}
If $\K$ is a local field
in the situation
of Lemma~\ref{cpinnormed},
it suffices to consider
\emph{ultrametric} continuous seminorms in~(b)
(as the proof shows).
\end{rem}
Recall that a topological vector space
$E$ over a topological field~$\K$
is called \emph{sequentially complete}
if every Cauchy sequence in $E$ is convergent.
We say that $E$ is \emph{Mackey complete}
if every Mackey-Cauchy sequence in~$E$
is convergent.
Here, a sequence $(x_n)_{n\in \N}$ in~$E$
is called a \emph{Mackey-Cauchy sequence}
if there exists a bounded subset $B\sub E$
and elements $\mu_{n,m}\in \K$
such that $x_n-x_m\in \mu_{n,m}B$ for all $n,m\in \N$
and $\mu_{n,m}\to 0$ in~$\K$ as both
$n,m\to\infty$.\\[2.5mm]
Note that every Mackey-Cauchy sequence also is
a Cauchy sequence; hence every sequentially complete
topological $\K$-vector space is
Mackey complete.
In the real locally convex case,
Mackey completeness is
a (particularly weak) standard
completeness property,
which is of great
usefulness for infinite-dimensional
calculus (see \cite[notably \S2]{KaM}
for an in-depth discussion).
%
%
%
\begin{thm}\label{holvsweak}
Let $\K\not=\C$ be a locally compact field,
$E$ and~$F$ be topological
$\K$-vector spaces, $f\colon U\to F$
be a map on an open set $U\sub E$,
$k\in \N_0\cup\{\infty\}$
and $\sigma\in \;]0,1]$.
If $E$ is metrizable
and $F$ is both Mackey complete
and locally convex,
then $f$ is $C^{k,\sigma}$
if and only if $f$ is weakly $C^{k,\sigma}$.
\end{thm}
\begin{proof}
The other implication being trivial,
we only need to show that if $f$ is weakly
$C^{k,\sigma}$, then $f$ is $C^{k,\sigma}$.
As a consequence of Theorem~\ref{dfc},
$f$ will be $C^{k,\sigma}$
if we can show that $g:=f\circ\gamma\colon \K^\ell\to F$
is $C^{k,\sigma}$ for each $\ell\in \N$
and each smooth map $\gamma \colon \K^\ell \to U$.
Note that $g$ is weakly $C^{k,\sigma}$
since so is~$f$ (see Remark~\ref{simplwea2}).
Hence, after replacing~$f$ with~$g$,
we may assume that $U=E=\K^\ell$ for some
$\ell\in \N$.
We may assume that $k\in \N_0$;
the proof is by induction on~$k$.\\[2.5mm]
If $k=0$ and $f\colon E=\K^\ell\to F$ is weakly $C^{0,\sigma}$,
let $x\in E$ and
$K \sub E$ be a compact neighbourhood
of~$x$. Then $f|_K$ is $C^{0,\sigma}$
by Lemma~\ref{cpinnormed}.
Hence $f$ is $C^{0,\sigma}$ locally
and hence $f$ is $C^{0,\sigma}$, by
Lemma~\ref{hoelloc}.\\[2.5mm]
\emph{Induction step}.
If $k\geq 1$ and $f\colon  E=\K^\ell\to F$
is weakly $C^{k,\sigma}$,
given $x,y\in E$
choose a sequence $(t_n)_{n\in \N}$ of pairwise distinct elements
in $\wb{B}^\K_1(0) \setminus \{0\}$
such that $t_n\to 0$. Set
\[
B\; :=\; \Big\{
\frac{f^{]1[}(x,y,t_m)
-f^{]1[}(x,y,t_n)}{\beta(|t_m-t_n|)} \colon n,m\in \N\Big \}\, ,
\]
where $\beta\colon \;]0,\infty[\, \to \K^\times$
is as in the proof of Lemma~\ref{cpinnormed}
(as well as $\rho$ used to define~$\beta$).
Then $\lambda(B)\sub \K$
is bounded for each $\lambda\in F'$
and hence
$B$ is bounded
(by \cite[Theorem~3.18]{Rud}, resp.,
\cite[Theorem~4.21]{Thi}).
In fact, since $\lambda\circ f$ is
$C^{1,\sigma}$, it follows
that $(\lambda\circ f)^{[1]}$
is
$C^{0,\sigma}$.
Applying now Lemma~\ref{cpinnormed}
to the restriction of
$(\lambda\circ f)^{[1]}$
to the compact set $\{x\}\times\{y\}\times \wb{B}^\K_1(0)$,
we find $C\in [0,\infty[$
such that
\[
|(\lambda\circ f)^{[1]}(x,y, t)
- (\lambda\circ f)^{[1]}(x,y, s)|\;\leq\;
C|t-s|^\sigma\quad
\mbox{for all $s,t\in \wb{B}^\K_1(0)$.}
\]
Repeating the calculation
in~(\ref{mayb}),
we find that $\sup\, |\lambda(B)|\leq C$.
Hence $B$ is indeed bounded.\\[2.5mm]
Since $f^{]1[}(x,y,t_m)
-f^{]1[}(x,y,t_n) \in \beta(|t_m-t_n|) B$,
where $B$ is bounded and
$\beta(|t_m-t_n|)\to 0$ as both $n,m\to\infty$,
we deduce that $(f^{]1[}(x,y,t_n))_{n\in \N}$
is a Mackey-Cauchy sequence in~$F$ and thus convergent;
we let $g(x,y,0)$ be its limit.
Then $\lambda(g(x,y,0))=\lim_{n\to\infty}(\lambda\circ f)^{]1[}(x,y,t_n)=
(\lambda\circ f)^{[1]}(x,y,0)$
for each $\lambda$. Furthermore, trivially
$\lambda(g(x,y,t))=(\lambda\circ f)^{[1]}(x,y,t)$
for $g(x,y,t):=f^{]1[}(x,y,t)=t^{-1}(f(x+ty)-f(x))$
if $(x,y,t)\in E\times E\times \K^\times$.
Thus $\lambda\circ g=(\lambda\circ f)^{[1]}$
is $C^{k-1,\sigma}$ for each $\lambda$, whence
$g$ is $C^{k-1,\sigma}$,
by induction. Hence~$f$ is~$C^{1,\sigma}$,
with $f^{[1]}=g$ a $C^{k-1,\sigma}$-map.
As a consequence,
$f$ is $C^{k,\sigma}$.
\end{proof}
\appendix
\section{Details for Section~\ref{secnot}}\label{appbore}
%
%
In this appendix, proofs are provided
for the lemmas of Section~\ref{secnot}.\\[2.5mm]
{\bf Proof of Lemma~\ref{presymm}.}
The assertions are obvious from
our definitions of $U^{>\alpha<}$
and $f^{>\alpha<}$.\vspace{2.5mm}\Punkt

\noindent
{\bf Proof of Lemma~\ref{interpr}.}
The $1$-dimensional case of this
lemma is well known (see \cite[Exercise~29.A]{Sch}).
Having done this exercise (or not),
the reader should not have
difficulties to work out the details
of the following sketch:
We start with formula (\ref{badbd}) for
$f^{>\alpha<}(x)$, and split the sum
$\sum_{j_i=0}^{\alpha_i}$
occurring there
into the sum $\sum_{j_i=1}^{\beta_i}$,
plus the two remaining summands with $j_i=0$
and $j_i=\alpha_i$, respectively.
In the sum $\sum_{j_i=1}^{\beta_i}$, rewrite
the factor $\frac{1}{x^{(i)}_{j_i}-x^{(i)}_0}\cdot
\frac{1}{x^{(i)}_{j_i}-x^{(i)}_{\alpha_i}}$\vspace{-.9mm}
of the product involved in the summands as
\[
\frac{1}{x^{(i)}_0-x^{(i)}_{\alpha_i}}\cdot
\left(
\frac{1}{x^{(i)}_{j_i}-x^{(i)}_0}-
\frac{1}{x^{(i)}_{j_i}-x^{(i)}_{\alpha_i}}\right).
\]
Finally, take $\frac{1}{x^{(i)}_0-x^{(i)}_{\alpha_i}}$
out of the sum
and combine the summands to
\begin{eqnarray*}
\lefteqn{f^{>\beta<}(x^{(1)},\ldots, x^{(i-1)},
x^{(i)}_0, \ldots, x^{(i)}_{\beta_i},
x^{(i+1)},\ldots, x^{(d)})}\qquad \\
& & -
f^{>\beta<}(x^{(1)},\ldots, x^{(i-1)},
x^{(i)}_1, \ldots, x^{(i)}_{\alpha_i},
x^{(i+1)},\ldots, x^{(d)})\, .
\end{eqnarray*}
After a reordering of the arguments in the second
term with the help of Lemma~\ref{presymm},
we obtain~(\ref{schlmm}).\,\vspace{2.5mm}\Punkt

\noindent
{\bf Proof of Lemma~\ref{partsymm}.}
The validity of (\ref{nweq})
is clear from the definition of $U^{<\alpha>}$.
Since $U^{>\alpha<}$ is dense in $U^{<\alpha>}$,
if suffices to check (\ref{eqsmm})
for $x\in U^{>\alpha<}$.
But then~(\ref{eqsmm}) holds by
Lemma~\ref{presymm}.\vspace{2.5mm}\Punkt

\noindent
{\bf Proof of Lemma~\ref{extmore}.}
We let $W$ be the set of all
$x\in U^{<\alpha>}$
such that $x^{(i)}_0\not=x^{(i)}_{\alpha_i}$,
and define $h(x)$ by~(\ref{schlmm2})
for $x\in W$.
Then  both $f^{<\alpha>}|_W$ and
$h\colon W\to F$
are continuous and coincide with $f^{>\alpha<}$
on $U^{>\alpha<}$,
by Lemma~\ref{interpr}.
Since $U^{>\alpha<}$ is dense in~$W$,
it follows that $f^{<\alpha>}|_W=h$.\vspace{2.5mm}\Punkt

\noindent
{\bf Proof of Lemma~\ref{bdviagau}.}
If $B\sub E$ is bounded and $q$
is a gauge, then there exists $t\in \K^\times$
such that $tB\sub B^q_1(0)$.
Thus $|t|\cdot \|x\|_q=\|tx\|_q<1$ for all
$x\in B$ and thus $\sup q(B)\leq |t|^{-1}$.
Conversely, suppose that
$q(B)$ is bounded for each gauge~$q$.
If $U\sub E$ is a $0$-neighbourhood,
there exists a gauge $q$ on~$E$
such that $B^q_1(0)\sub U$.
Choose $r>0$ such that $r\sup q(B)<1$.
Then $B^{\K}_r(0) \cdot B\sub B_1^q(0)\sub U$,
showing that~$B$ is bounded.\,\vspace{2.5mm}\Punkt

\noindent
{\bf Proof of Lemma~\ref{baseHglob}.}
(a)
Given $x_0\in U$, a gauge $q$ on~$F$,
and $\varepsilon>0$,
we choose a gauge $p$ on~$E$ and a neighbourhood $W\sub U$ of~$x_0$
such that\linebreak
$\|f(y)-f(x)\|_q\leq (\|y-x\|_p)^\sigma$
for all $x,y\in W$.
Define $\delta:=\ve^\frac{1}{\sigma}$.
Then\linebreak
$f(W \cap B^p_\delta(x_0))\sub B_{\delta^\sigma}^q(f(x_0))
=B_\ve^q(f(x_0))$, as
$\|f(y)-f(x_0)\|_q\leq (\|y-x_0\|_p)^\sigma$
for all $y\in W$.
Hence $f$ is continuous at $x_0$
(see \cite[Lemma~1.27\,(b)]{IM2}).

(b) Since $f$ is $C^{0,\sigma}$, given
$x_0\in U$ and a gauge~$q$ on~$F$,
we find a gauge $p$ on~$E$ and a neighbourhood $W \sub U$ of~$x_0$
such that
$\|f(y)-f(x)\|_q\leq $ $(\|y-x\|_p)^\sigma$
for all $x,y\in W$.
Let~$s$ be a gauge on~$E$ such that
$p(x+y)\leq s(x)+s(y)$ for $x,y\in E$.
After replacing $W$ with $W \cap B^s_{1/2}(x_0)$,
we may assume that $\|y-x\|_p<1$
for all $x,y\in W$. Then
$\|f(y)-f(x)\|_q\leq (\|y-x\|_p)^\sigma
= (\|y-x\|_p)^{\sigma-\tau}
(\|y-x\|_p)^\tau
\leq (\|y-x\|_p)^\tau$
for all $x,y\in W$,
whence~$f$ is~$C^{0,\tau}$. 

(c) Let $x_0\in U$. Given a gauge $q$ on~$H$,
there exists a gauge~$p$ on~$F$ and a neighbourhood
$R\sub V$ of $f(x_0)$ such that
$\|g(y)-g(x)\|_q\leq (\|y-x\|_p)^\tau$
for all $x,y\in R$.
There exists a gauge~$s$ on~$E$ and a
neighbourhood $S\sub f^{-1}(R)$ of~$x_0$
such that $\|f(y)-f(x)\|_p\leq(\|y-x\|_s)^\sigma$
for all $x,y\in S$.
Then $\|g(f(y))-g(f(x))\|_q\leq
(\|f(y)-f(x)\|_p)^\tau
\leq (\|y-x\|_s)^{\sigma\cdot \tau}$
for all $x,y\in S$.

(d) See \cite[Lemma~2.5\,(c)]{IM2}.\vspace{2mm}\Punkt

\noindent
{\bf Proof of Lemma~\ref{densehoel}.}
If $f$ is not $C^{0,\sigma}$,
then there exists $x_0\in U$
and a gauge~$q_0$ on~$F$ such that,
for each neighbourhood $V\sub U$
of~$x_0$ and gauge~$p$ on~$E$,
there are $x,y\in V$ such that
$\|f(y)-f(x)\|_{q_0}>
(\|y-x\|_p)^\sigma$.
Let~$q$
be a gauge on~$F$
such that $q_0(u+v)\leq
q(u)+q(v)$
for all $u,v\in F$.
After replacing~$q$ with a larger
gauge, we may assume that~$q$ is upper semicontinuous
(cf.\ Remark~\ref{minkow} and Lemma~\ref{findfund}).
We now verify that $x_0$
and~$q$ have the desired properties.
To see this, let $V\sub U$ be a
neighbourhood of~$x_0$
and~$p_0$ be a gauge on~$E$.
Let $p\geq p_0$ be an upper semicontinuous
gauge.
Then there are $x,y\in V$
such that
$\ve:=
\|f(y)-f(x)\|_{q_0}-(\|y-x\|_p)^\sigma>0$.
Choose $r>\|y-x\|_p$
such that $r^\sigma\leq (\|y-x\|_p)^\sigma+\frac{\ve}{2}$.
Since $B^p_r(0)$ and $B^q_{\ve/2}(0)$
are open
and the relevant maps are continuous,
we find $x',y'\in V\cap D$
such that $\|y'-x'\|_p<r$
and $\|f(y)-f(x)-f(y')+f(x')\|_q
<\frac{\ve}{2}$.
Using the fake triangle inequality,
we now obtain
\begin{eqnarray*}
\|f(y')-f(x')\|_q
& \geq &
\|f(y)-f(x)\|_{q_0}-\|f(y)-f(x)-f(y')+f(x')\|_q\\
& > &
\|f(y)-f(x)\|_{q_0}-\frac{\ve}{2}
\; =\;
(\|y-x\|_p)^\sigma +\frac{\ve}{2}\\
& \geq &
(\|y'-x'\|_p)^\sigma\;\geq \; (\|y'-x'\|_{p_0})^\sigma\,,
\end{eqnarray*}
as desired.\vspace{2.5mm}\Punkt

\noindent
{\bf Proof of Lemma~\ref{chainlin}.}
If $k=0$, then $\lambda\circ f$ is $C^{0,\sigma}$
by Lemma~\ref{baseHglob}\,(c),\linebreak
exploiting that $\lambda$, being continuous
linear, is $C^\infty_{BGN}$ and hence Lipschitz continuous.
If $f$ is $C^{k+1,\sigma}_{BGN}$,
then $\lambda\circ f$ is $C^{k,\sigma}_{BGN}$
by induction, with $(\lambda\circ f)^{[k]}
=\lambda\circ f^{[k]}$.
Furthermore, $\lambda\circ f$ is $C^{k+1}_{BGN}$,
by~{\bf\ref{chainr}}.
Now $(\lambda\circ f)^{[k+1]}=
((\lambda\circ f)^{[k]})^{[1]}=\lambda^{[1]}\circ
\wh{T}(f^{[k]})=
\lambda \circ (f^{[k]})^{[1]}=\lambda\circ f^{[k+1]}$
is $C^{0,\sigma}$,
using~{\bf\ref{chainr}}
for the second equality,
{\bf\ref{exalin}} for the third
(cf.\ also \cite[Remark~1.7]{IM2}).
Hence $\lambda \circ f$ is $C^{k+1,\sigma}_{BGN}$
with $(\lambda\circ f)^{[k+1]}=\lambda\circ
f^{[k+1]}$ of the desired form.\vspace{2.5mm}\Punkt

\noindent
{\bf Proof of Lemma~\ref{prodhoel}.}
If $f$ is $C^{k,\sigma}_{BGN}$,
then also $f_i=\pr_i\circ f$
is $C^{k,\sigma}_{BGN}$, since $\pr_i$ is continuous linear
(Lemma~\ref{chainlin}).
Conversely, assume that each component~$f_i$
of $f\colon U\to \prod_{i\in I}F_i=F$
is $C^{k,\sigma}_{BGN}$.
We proceed by induction.\\[2.5mm]
\emph{The case $k=0$.}
Given a gauge~$q$ on $F$,
there exists a finite subset $J\sub I$
and balanced, open $0$-neighbourhoods
$W_j\sub F_j$ for $j\in J$
such that $W:=\bigcap_{j\in J}\pr_j^{-1}(W_j)\sub
B^q_1(0)$.
We let $s:=\mu_W\colon F\to [0,\infty[$ be the Minkowski\linebreak
functional
of~$W$ (see Remark~\ref{minkow}),
and $s_j\colon F_j \to [0,\infty[$
be the Minkowski functional
of $W_j$, for $j\in J$.
Then $s(x)=\max\{s_j(x_j)\colon j\in J\}$
holds for each $x=(x_i)_{i\in I}\in F$.
Furthermore, $q(x) \leq s(x)$.
In fact, given $x\in F$ and $t\in \K^\times$
such that
$x\in t W$, we have
$q(x)=q(t(x/t))=|t|\cdot q(x/t)\leq |t|$
(using that $W\sub B^q_1(0)$).
Letting $|t|\to\mu_W(x)$, we see that
$q(x)\leq \mu_W(x)=s(x)$.
For each $j\in J$,
there is a gauge $p_j$ on~$E$
and a neighbourhood $V_j$ of~$x_0$ in~$U$
such that $\|f_j(y)-f_j(x)\|_{s_j}\leq
(\|y-x\|_{p_j})^\sigma$
for all $x,y\in V_j$.
Set $V:=\bigcap_{j\in J} V_j$
and $p(x):=\max\{p_j(x)\colon j\in J\}$
for $x\in E$.
Then~$p$ is a gauge on~$E$ such that
$\|f(y)-f(x)\|_q\leq \|f(y)-f(x)\|_s
=\max\{\|f_j(y)-f_j(x)\|_{s_j}\colon j\in J\}
\leq (\|y-x\|_p)^\sigma$
for all $x,y\in V$.\\[2.5mm]
\emph{Induction step.}
Assume that each component~$f_i$
is $C^{k,\sigma}_{BGN}$, where $k\geq 1$.
Then $f$ is $C^1_{BGN}$, with
$f^{[1]}=(f_i^{[1]})_{i\in I}$
(cf.\ \cite[Lemma~10.2]{BGN}).
The components $f_i^{[1]}$ of this map
are $C^{k-1,\sigma}_{BGN}$,
whence $f^{[1]}$ is $C^{k-1,\sigma}_{BGN}$,
by induction.
Hence $f$ is $C^{k,\sigma}_{BGN}$.\vspace{2.5mm}\Punkt

\noindent
{\bf Proof of Lemma~\ref{chainLip}.}
We proceed by induction on $k\in \N_0$.
If $k=0$, then Lemma~\ref{chainLip}
is a special case of Lemma~\ref{baseHglob}\,(c).\\[2.5mm]
\emph{Induction step}:
If $f$ is $C^{k,\sigma}_{BGN}$
and $g$ is $C^{k,\tau}_{BGN}$ with $k\geq 1$,
then $g\circ f$ is $C^k_{BGN}$
and $(g\circ f)^{[1]}=g^{[1]}\circ \wh{T}f$
with $\wh{T}f\colon U^{[1]}\to V^{[1]}\sub F\times F\times \K$,
$\wh{T}f(x,y,t):= (f(x),f^{[1]}(x,y,t), t)$
(see {\bf\ref{chainr}}).
Here the second component of $\wh{T}f$ is $C^{k-1,\sigma}_{BGN}$;
the final component is continuous linear and hence
$C^{\infty,\sigma}_{BGN}$ (since $\sigma\leq 1$); and the first
component is a composition of the $C^{k,\sigma}_{BGN}$-map
$f$ and (a restriction of) the continuous
linear (and hence $C^{\infty,\sigma}_{BGN}$-)
mapping $E\times E\times \K\to E$,\linebreak
$(x,y,t)\mto x$,
whence also the first component is $C^{k-1,\sigma}_{BGN}$,
by the case $k-1$
(valid by induction).
Now Lemma~\ref{prodhoel}
shows that $\wh{T}f$
is $C^{k-1,\sigma}_{BGN}$,
and thus
$(g\circ f)^{[1]}=g^{[1]}\circ \wh{T}f$
is $C^{k-1,\sigma\cdot \tau}_{BGN}$,
by induction.
Hence $g\circ f$
is $C^{k,\sigma\cdot \tau}_{BGN}$.\vspace{2.5mm}\Punkt

\noindent
{\bf Proof of Lemma~\ref{hoelloc}.}
We may assume that $k\in \N_0$;
the proof is by induction.
If $k=0$ and $f|_{U_i}$ is $C^{0,\sigma}$
for each $i\in I$, then $f$ is $C^{0,\sigma}$,
as is obvious from the definition.
Induction step: If $f|_{U_i}$ is $C^{k+1,\sigma}_{BGN}$,
then $f$ is $C^{k,\sigma}_{BGN}$ by induction,
and furthermore $f$ is~$C^1_{BGN}$ (by Lemma~\ref{Crlocal}).
The sets $U_i^{[1]}$ together with
$U^{]1[}$ form an open cover
for $U^{[1]}$, and
$f^{[1]}|_{U_i^{[1]}}=(f|_{U_i})^{[1]}$
is $C^{k,\sigma}_{BGN}$ for each $i\in I$.
For $(x,y,t)\in U^{]1[}$, we have
$f^{[1]}(x,y,t)=\frac{f(x+ty)-f(x)}{t}$;
since $f$ is $C^{k,\sigma}_{BGN}$, we deduce with
Lemma~\ref{chainLip}
from the preceding formula
that $f^{[1]}|_{U^{]1[}}$
is $C^{k,\sigma}_{BGN}$.
Applying the inductive hypothesis,
we see that $f^{[1]}$ is $C^{k,\sigma}_{BGN}$.
Hence $f$ is $C^{k+1,\sigma}_{BGN}$.\,\Punkt
\section{Topologies on function spaces}\label{apptop}
%
%
Let $\K$ be a topological field, $F$ be a topological
$\K$-vector space, $d\in \N$,
$k\in \N_0 \cup \{\infty\}$
and $U\sub \K^d$ be a subset
which is open or of the form
$U=U_1 \times\cdots\times U_d$
for certain subsets $U_1,\ldots, U_d\sub \K$
with dense interior.
Since $C^k_{BGN}$-maps
and $C^k_{SDS}$-maps $U\to F$
coincide (by Theorem~\ref{mainthm}), as before we shall simply
refer to them as $C^k$-maps.
Now both the definition
of $C^k_{BGN}$-maps and the definition
of $C^k_{SDS}$-maps suggest a natural
definition of a vector topology
on the space $C^k(U,F)$ of all
$C^k$-maps $U\to F$.
The first of these has been used in~\cite{ZOO};
the second one has been used in~\cite{Sch}
and further works by Schikhof,
in the case of functions of a single variable.
In this section,
we show that the two topologies
coincide.
%
%
\begin{defn}\label{deftps}
We write $C^k(U,F)_{BGN}$
for the space $C^k(U,F)$, equipped with the
initial topology with respect to the
maps
%
\begin{equation}\label{dftau}
\tau_j\colon C^k(U,F)\to C(U^{[j]},F)\, ,
\quad f\mto f^{[j]}
\end{equation}
for $j\in \N_0$ such that $j\leq k$,
where $C(U^{[j]},F)$ is equipped with the
compact-open topology (which coincides
with the topology of uniform convergence on compact
sets).
We write $C^k(U,F)_{SDS}$
for $C^k(U,F)$, equipped with the
initial topology with respect to the
maps
%
\begin{equation}\label{dfkapp}
\kappa_\alpha \colon C^k(U,F)\to C(U^{<\alpha>},F)\, ,
\quad f\mto f^{<\alpha>}
\end{equation}
for $\alpha \in \N_0^d$ such that $|\alpha| \leq k$,
where $C(U^{<\alpha>},F)$ is equipped with the
compact-open topology.
\end{defn}
It is clear that both
$C^k(U,F)_{BGN}$ and $C^k(U,F)_{SDS}$
are (Hausdorff) topological vector spaces
over~$\K$.
Information concerning basic
properties of $C^k(U,F)_{BGN}$
(like completeness, metrizability,
and local convexity) can be found
in \cite[Proposition~4.19]{ZOO}, for~$U$ an open
subset of a topological $\K$-vector
space (or more generally
a $C^k$-manifold modelled
on a topological $\K$-vector space).
Our goal is the following theorem.
%
\begin{thm}\label{thmsametp}
$C^k(U,F)_{BGN}=C^k(U,F)_{SDS}$
as a topological $\K$-vector space.
\end{thm}
It is essential for the proof of
Theorem~\ref{thmsametp} that the maps~$f^{[j]}$
and~$f^{<\alpha>}$ can be expressed
in terms of each other.
We do not need explicit formulas;
the information provided by
the next two lemmas is sufficient
for our purposes, and easy to work with.
%
%
\begin{numba}\label{prplemA}
To prepare the proof of following lemma,
consider a continuous affine-linear map
$\theta\colon X\to Y$
between topological $\K$-vector spaces~$X$ and~$Y$,
say $\theta(x)=\lambda(x)+y_0$
with $y_0\in Y$ and a continuous linear
map $\lambda\colon X\to Y$.
Then $\theta^{[1]}(x,y,t)=\lambda(y)$
(cf.\ {\bf\ref{exalin}}), entailing that
$\wh{T}\theta\colon X^{[1]}\to Y^{[1]}$ (as in {\bf\ref{chainr}})
is given by
$(\wh{T}\theta)(x,y,t)=(\theta(x),\lambda(y),t)$
and hence is a continuous affine-linear map.
\end{numba}
%
%
\begin{la}\label{alphviaquot}
For each $\alpha\in \N_0^d$
such that $j:=|\alpha|\leq k$,
there exists an
affine-linear map
$\theta_\alpha\colon (\K^d)^{<\alpha>}\to (\K^d)^{[j]}$
such that $\theta_\alpha (U^{<\alpha>})\sub U^{[j]}$
and
\[
f^{<\alpha>}\; =\;\, f^{[j]}\circ \, \theta_\alpha|_{U^{<\alpha>}}
\]
for each $C^k$-map $f\colon U\to F$.
\end{la}
\begin{proof}
The proof is by induction on $j:=|\alpha|$.
If $j=0$, then $U^{[j]}=U^{<\alpha>}=U$
and $f^{<\alpha>}=f=f^{[j]}=f^{[j]}\circ
\theta_\alpha|_{U^{<\alpha>}}$
with $\theta_\alpha:=\id\colon \K^d\to \K^d$, $x\mto x$
linear.\\[2.5mm]
\emph{Induction step.}
If $j:=|\alpha|\geq 1$, then $\alpha=\beta+e_i$
for some $\beta\in \N_0^d$ and some
$i\in \{1,\ldots, d\}$.
Define $h \colon (\K^d)^{<\alpha>}\to ((\K^d)^{<\beta>})^{[1]}$
via
\[
h(x) :=
(x^{(1)},\ldots, x^{(i-1)},
x^{(i)}_0,\ldots, x^{(i)}_{\beta_i},
x^{(i+1)},\ldots,x^{(d)}, e_{\alpha_1+\cdots+\alpha_{i-1}+i}
, x^{(i)}_{\alpha_i}
-x^{(i)}_0)
\]
for $x\in (\K^d)^{<\alpha>}$,
where $e_{\alpha_1+\cdots+\alpha_{i-1}+i}$ is the
indicated unit vector in~$(\K^d)^{<\beta>}$.
Given $x\in U^{>\alpha<}$
(resp., $x\in U^{<\alpha>}$),
it is clear from the definitions
that
$h(x)
\in (U^{<\beta>})^{]1[}$
(resp.,
$h(x)\in (U^{<\beta>})^{[1]}$),
and
%
\begin{equation}\label{pa1eqn}
f^{<\alpha>}(x)\; =\; (f^{<\beta>})^{[1]}(h(x))\quad
\mbox{for all $x\in U^{>\alpha <}$.}
\end{equation}
Hence
$f^{<\alpha>}(x)=(f^{<\beta>})^{[1]}(h(x))$
for all $x\in U^{<\alpha>}$,
by continuity of the functions involved and density
of~$U^{> \alpha <}$ in~$U^{<\alpha>}$.
Here $f^{<\beta>}=f^{[j-1]}\circ \theta_\beta|_{U^{<\beta>}}$
by induction, with $\theta_\beta\colon
(\K^d)^{<\beta>} \to
(\K^d)^{[j-1]}$ an affine-linear map.
Then
%
\begin{equation}\label{pa2eqn}
(f^{<\beta>})^{[1]}
\; =\; f^{[j]}\circ (\wh{T}\theta_\beta)|_{(U^{<\beta>})^{[1]}}\, ,
\end{equation}
where $\wh{T}\theta_\beta\colon
((\K^d)^{<\beta>})^{[1]}\to
(\K^d)^{[j]}$
is affine-linear by {\bf\ref{prplemA}}.
Combining~(\ref{pa1eqn}) and~(\ref{pa2eqn})
gives $f^{<\alpha>}=f^{[j]}\circ \theta_\alpha|_{U^{<\alpha>}}$
with $\theta_\alpha:=(\wh{T}\theta_\beta)\circ h\colon
(\K^d)^{<\alpha>}\to (\K^d)^{[j]}$
an affine-linear map.
\end{proof}
%
%
%
\begin{numba}\label{prplemB}
We recall that a function
$p\colon \K^d\to F$ to a topological $\K$-vector space~$F$
is called a \emph{polynomial function}
if there exists $N\in \N_0$
and elements $a_\alpha\in F$ for each multi-index
$\alpha\in \N_0^d$ with $|\alpha|\leq N$,
such that
\[
f(x_1,\ldots, x_d)\;=\;\sum_{|\alpha|\leq N}a_\alpha\,
x_1^{\alpha_1}\cdots x_d^{\alpha_d}\quad
\mbox{for all $x_1,\ldots, x_d\in \K$.}
\]
It is easy to see that each polynomial function
$p\colon \K^d\to F$ is~$C^\infty$,
and that~$p^{[1]}$ and $\wh{T}p\colon (\K^d)^{[1]}\to F^{[1]}$
are polynomial functions
if so is~$p$.
\end{numba}
%
%
\begin{la}\label{quotviaalph}
Let $U=U_1\times\cdots\times U_d$
for certain subsets $U_1,\ldots, U_d\sub \K$
with dense interior.
Let $j\in \N_0$ such that $j\leq k$.
Then there exists $\ell\in \N$,
polynomial functions $p_1,\ldots, p_\ell\colon
(\K^d)^{[j]}\! \to\! \K$, multi-indices
$\alpha^{(1)},\ldots,\alpha^{(\ell)}\in \N_0^d$
and polynomial functions
$q_i \colon (\K^d)^{[j]}\to (\K^d)^{<\alpha^{(i)}>}$
for $i\in \{1,\ldots, \ell\}$
such that $q_i(U^{[j]})\sub U^{<\alpha^{(i)}>}$
for each~$i$ and
%
\begin{equation}\label{exprdff}
f^{[j]}(x) \;=\; \sum_{i=1}^\ell p_i(x)\cdot
f^{<\alpha^{(i)}>}(q_i(x))
\end{equation}
for each $f\in C^k(U,F)$ and each $x\in U^{[j]}$,
i.e.,
\begin{equation}\label{exprdff2}
f^{[j]} \;=\;
\sum_{i=1}^\ell p_i|_{U^{[j]}}\cdot
(f^{<\alpha^{(i)}>}\circ q_i|_{U^{[j]}})\, .
\end{equation}
\end{la}
\begin{proof}
The proof is by induction on~$j$.
If $j=0$, then $U^{<0>}=U=U^{[0]}$
and $f^{<0>}=f=f^{[0]}$,
whence $f^{[0]}=(p|_{U^{[0]}})\cdot (f^{<0>}\circ
q|_{U^{[0]}})$ with $p(x):=1$ and
$q:=\id\colon \K^d\to \K^d$.\\[2.5mm]
Now suppose that~$f^{[j]}$ is of the form~(\ref{exprdff2})
for some $j\in \N_0$ such that $j+1\leq k$.
Using the Product Rule
\cite[\S\,3.3]{BGN},
we deduce from~(\ref{exprdff2}) that
%
%
\begin{eqnarray}
f^{[j+1]}(x,y,t) & =\ &
\sum_{i=1}^\ell\, \big(
(p_i)^{[1]}|_{U^{[j+1]}}(x,y,t)\cdot f^{<\alpha^{(i)}>}(q_i(x))\notag\\
& & \;\;\; + \, p_i(x) \cdot
(f^{<\alpha^{(i)}>}\circ q_i|_{U^{[j]}})^{[1]}(x,y,t)\notag \\
& & \;\;\;
+ \, t \cdot (p_i)^{[1]}(x,y,t)\cdot
(f^{<\alpha^{(i)}>}\circ q_i|_{U^{[j]}})^{[1]}(x,y,t)\big)
\label{goodpieces}
\end{eqnarray}
for all $(x,y,t)\in U^{[j+1]}\sub U^{[j]}\times
E^{[j]}\times \K$.
By~{\bf\ref{chainr}}, we have
%
\begin{equation}\label{goop2}
(f^{<\alpha^{(i)}>}\circ q_i|_{U^{[j]}})^{[1]}
\;=\;
(f^{<\alpha^{(i)}>})^{[1]}
\circ (\wh{T}q_i)|_{U^{[j+1]}}
\end{equation}
here. By {\bf\ref{prplemB}},
$(p_i)^{[1]}$
and $\wh{T}q_i\colon (\K^d)^{[j+1]}\to
((\K^d)^{<\alpha^{(i)}>})^{[1]}$
are polynomial functions.
In view of~(\ref{goodpieces})
and~(\ref{goop2}), and exploiting the
fact that compositions of polynomial functions
are polynomial,
$f^{[j+1]}$ will be of the desired
form (analogous to~(\ref{exprdff2}))
if we can show that
$(f^{<\alpha^{(i)}>})^{[1]}(x)$
for $x\in (U^{<\alpha^{(i)}>})^{[1]}$
is a sum of terms of the form
\[
Q_m(x) f^{<\alpha^{(i)}+\wb{e_m}>}(P_m(x))
\]
for $m\in \{1,\ldots, d+|\alpha^{(i)}|\}$
and polynomial functions
$Q_m\colon ((\K^d)^{<\alpha^{(i)}>})^{[1]}\to\K$
and $P_m\colon ((\K^d)^{<\alpha^{(i)}>})^{[1]}\to
(\K^d)^{<\alpha^{(i)}+\wb{e_m}>}$
with
$P_m(U^{<\alpha^{(i)}>})^{[1]})\sub
U^{<\alpha^{(i)}+\wb{e_m}>}$,
where $e_m\in \K^{d+|\alpha^{(i)}|}$
is a standard unit vector
and $\wb{e_m}\in \K^d$
as in Remark~\ref{remeazy}.
For simplicity of notation,
write $g:=f^{<\alpha^{(i)}>}$
and $N:=d+|\alpha^{(i)}|$.
Since $U^{<\alpha^{(i)}>}=(U_1)^{\alpha^{(i)}_1}\times\cdots
\times (U_d)^{\alpha^{(d)}_d}\sub \K^N$
is a product of subsets of~$\K$,
we have
\begin{eqnarray*}
\lefteqn{g^{[1]}(x, y,t)}\\
&= & \sum_{m=1}^N y_m\cdot
g^{<e_m>}(x_1+ty_1,\ldots,
x_{m-1}+ty_{m-1};x_m,x_m+ty_m;
x_{m+1},\ldots, x_N)
\end{eqnarray*}
for all $(x,y,t) \!\in\!
(U^{<\alpha^{(i)}>})^{[1]}$
with $x\!=\! (x_1,\ldots, x_N)\! \in\! \K^N$,
$y\!=\! (y_1,\ldots, y_N)\!\in\! \K^N$
and $t\in\K$ (cf.\ proof of Theorem~\ref{mainthm}).
Since
$g^{<e_m>}= (f^{<\alpha^{(i)}>})^{<e_m>}
=f^{<\alpha^{(i)}+ \wb{e_m}>}$
with notation as in Remark~\ref{remeazy},
we deduce that $(f^{<\alpha^{(i)}>})^{[1]}=g^{[1]}$
is a sum of terms of the desired form.
\end{proof}
%
%
%
\begin{numba}\label{contfctsp}
We recall:
If $X,Y,Z$ are Hausdorff topological spaces
and $g\colon X\to Y$
is a continuous map,
then
\[
C(Z,g)\colon C(Z,X)\to C(Z,Y)\, ,\quad
f\mto g\circ f
\]
and
\[
C(g,Z)\colon C(Y,Z)\to C(X,Z)\,,\quad
f\mto f\circ g
\]
are continuous maps, using the compact-open
topology on all function spaces
(see \cite[p.\,157, Assertions (1) and (2)]{Eng}).
If $X$ is a Hausdorff topological space
and~$F$ a topological $\K$-vector space,
then $C(X,F)$ is a topological $C(X,\K)$-module
with pointwise operations
(by continuity of $C(X,\alpha)$ and $C(X,\mu)$, 
where
$\mu\colon \K\times F\to F$, $\mu(t,x):=tx$
is scalar multiplication
and $\alpha \colon F\times F\to F$,
$\alpha(x,y):=x+y$).
This implies that the multiplication operator
\[
m_h\colon C(X,F)\to C(X,F)\,,\quad
f\mto h\cdot f
\]
is a continuous linear map for each $h\in C(X,\K)$.
\end{numba}
We also need to understand maps of the form $C^k(g,F)$.
%
%
\begin{la}\label{puback}
Let $d,e\in \N$,
$U\sub \K^d$ be an open subset
of the form $U=U_1\times \cdots\times U_d$
for certain open sets $U_1,\ldots,U_d\sub \K$,
$V\sub \K^e$ be open
and $g\colon U\to V$ be a $C^k_{SDS}$-map,
where $k\in \N_0\cup\{\infty\}$.
Then the map
\[
C^k(g,F)_{SDS}\colon
C^k(V,F)_{SDS}\to C^k(U,F)_{SDS}\, ,\quad
f\mto f\circ g
\]
is continuous and linear,
for each topological $\K$-vector space~$F$.
\end{la}
\begin{proof}
It is clear that
$C^k(g,F)_{SDS}$ is linear.
The map $C^k(g,F)_{SDS}$
will be continuous if we can show
that $\kappa_\alpha\circ C^k(g,F)_{SDS}$
is continuous for
each $\alpha\in\N_0^d$ such that
$|\alpha|\leq k$,
where $\kappa_\alpha$ is as in~(\ref{dfkapp}).
Set $j:=|\alpha|$.
Using Lemma~\ref{alphviaquot},
Lemma~\ref{quotviaalph}
and the notation introduced there,
we obtain
%
\begin{eqnarray}
(\kappa_\alpha\circ C^k(g,F)_{SDS})(f)
& =& (f\circ g)^{<\alpha>}
\;=\; (f\circ g)^{[j]}\circ \theta_\alpha|_{U^{<\alpha>}}\notag\\
& =& f^{[j]}\circ (\wt{T}^jg)\circ \theta_\alpha|_{U^{<\alpha>}}\notag\\
&=&
\sum_{i=1}^\ell f_i
\cdot
(f^{<\alpha^{(i)}>}\circ g_i)\label{strain}
\end{eqnarray}
with $f_i:=(p_i|_{U^{[j]}}) \circ
(\wt{T}^jg)\circ \theta_\alpha|_{U^{<\alpha>}}$
and $g_i:= (q_i|_{U^{[j]}}) \circ
(\wt{T}^jg)\circ \theta_\alpha|_{U^{<\alpha>}}$.
Using the continuous maps
$\eta_\beta\colon C^k(V,F)_{SDS}\to C(V^{<\beta>},F)$,
$f\mto f^{<\beta>}$,
the multiplication
operator $m_{f_i}\colon C(U^{<\alpha>},F)\to
C(U^{<\alpha>},F)$ by~$f_i$
for~$i$ in $\{1,\ldots, e\}$
and the maps
$C(g_i,F)\colon C(V^{<\alpha^{(i)}>},F)\to C(U^{<\alpha>},F)$
(which are continuous linear by {\bf\ref{contfctsp}}),
we can rewrite (\ref{strain}) in the form
\[
\kappa_\alpha\circ C^k(g,F)_{SDS}\;=\;
\sum_{i=1}^\ell \, m_{f_i}
\circ C(g_i,F)\circ \eta_{\alpha^{(i)}}\,.
\]
The right hand side being composed
of continuous maps, we deduce that
$\kappa_\alpha\circ C^k(g,F)_{SDS}$
is continuous.
\end{proof}
%
%
\begin{numba}\label{compincov}
Let $X$ and $Y$ be Hausdorff topological
spaces and $(U_i)_{i\in I}$ an open cover
of $X$. Then each compact set $K\sub X$
can be written as $K=\bigcup_{i\in \Phi}K_i$
for a certain finite subset $\Phi\sub I$
and compact sets $K_i\sub U_i$.
As a consequence,
the compact-open topology
on $C(X,Y)$ is the initial topology
with respect to the restriction maps
$C(X,Y)\to C(U_i,Y)$, $f\mto f|_{U_i}$
for $i\in I$,
where $C(U_i,Y)$ is equipped with the
compact-open topology.
\end{numba}
{\bf Proof of Theorem \ref{thmsametp}.}
The topology on $C^\infty(U,F)_{BGN}$
is initial with respect to the inclusion maps
$C^\infty(U,F)_{BGN}\to C^k(U,F)_{BGN}$
for $k\in \N_0$,
as is clear from the definitions.
Also, the topology on $C^\infty(U,F)_{SDS}$
is initial with respect to the inclusion maps
$C^\infty(U,F)_{SDS}\to C^k(U,F)_{SDS}$
for $k\in \N_0$.
As a consequence,
$C^\infty(U,F)_{BGN}=C^\infty(U,F)_{SDS}$
will hold if
$C^k(U,F)_{BGN}=C^k(U,F)_{SDS}$
for each $k\in \N_0$.
The proof is by induction on~$k$.\\[2.5mm]
If $k=0$, then both
$C^0(U,F)_{BGN}=C^0(U,F)_{SDS}=C(U,F)$,
equipped with the compact-open topology.\\[2.5mm]
Now let $k\in \N$
and suppose
that the topologies on spaces of $C^j$-maps
coincide, for all $j<k$.
To establish the equality $C^k(U,F)_{BGN}=C^k(U,F)_{SDS}$,
we prove that both of the following maps
are continuous:
\[
\Phi\colon C^k(U,F)_{BGN}\to C^k(U,F)_{SDS}\, ,
\quad f\mto f
\]
and
\[
\Psi\colon C^k(U,F)_{SDS}\to C^k(U,F)_{BGN}\, ,
\quad f\mto f\,.
\]
If we want to emphasize~$k$, we also
write $\Psi_k$ instead of~$\Psi$.\\[2.5mm]
\emph{Continuity of~$\Phi$:}
The topology on $C^k(U,F)_{SDS}$
being initial with respect to the maps
$\kappa_\alpha$ (from (\ref{dfkapp}))
for $\alpha\in \N_0^d$ such that
$|\alpha|\leq k$, the map $\Phi$ will
be continuous if we can show that $\kappa_\alpha\circ \Phi$
is continuous for each~$\alpha$.
Set $j:=|\alpha|$ and let~$\tau_j$ be as in~(\ref{dftau}).
By Lemma~\ref{alphviaquot},
we have $(\kappa_\alpha\circ \Phi)(f)=f^{<\alpha>}
=f^{[j]}\circ \theta_\alpha|_{U^{<\alpha>}}$
for each $f\in C^k(U,E)$.
Hence $\kappa_\alpha\circ \Phi=
C(\theta_\alpha|_{U^{<\alpha>}},F)\circ \tau_j$
is continuous by continuity of~$\tau_j$
and continuity of the map
$C(\theta_\alpha|_{U^{<\alpha>}},F)$
(see {\bf\ref{contfctsp}}).
Thus $\Phi$ is continuous.\\[2.5mm]
\emph{Continuity of~$\Psi$.}
We have to show that $\tau_j\circ \Psi$ is continuous
for each $j\in \N_0$ such that $j\leq k$.
For $j=0$, we have $\tau_0\circ \Psi=\kappa_0$,
which is continuous. We can therefore assume
now that $j\geq 1$.
Let us assume first that
$U=U_1\times\cdots\times U_d$,
where $U_1,\ldots,U_d\sub \K$
are subsets with dense interior.
For each $f\in C^k(U,F)$,
we have
$(\tau_j\circ \Psi)(f)=f^{[j]}
=\sum_{i=1}^\ell (p_i|_{U^{[j]}})
\cdot
(f^{<\alpha^{(i)}>}\circ q_i|_{U^{[j]}})$,
with notation as in Lemma~\ref{quotviaalph}.
Abbreviating $f_i:=(p_i|_{U^{[j]}})$,
we can rewrite the preceding formula as
%
\begin{equation}\label{showscts}
\tau_j\circ \Psi\; =\;
\sum_{i=1}^\ell \, m_{f_i}\circ
C(q_i|_{U^{[j]}},F)\circ
\kappa_{\alpha^{(i)}}\,,
\end{equation}
where $m_{f_i}\colon C(U^{[j]},F)\to
C(U^{[j]},F)$, $g\mto
f_i \cdot g$
is the multiplication operator
by~$f_i$.
Since all mappings on the right hand side of
(\ref{showscts}) are continuous
by~{\bf\ref{contfctsp}},
we deduce that
$\tau_j\circ \Psi$ is continuous
and hence also~$\Psi$.\\[2.5mm]
Now assume that $U\sub \K^d$ is an arbitrary
open set. Then there exists an open cover
$(U_i)_{i\in I}$ of~$U$
by sets of the form
$U_i=U_{i,1}\times \cdots\times U_{i,d}$,
where $U_{i,1}, \ldots, U_{i,d}\sub \K$
are open.
Now $(U^{[j-1]})^{]1[}$,
together with the sets $(U_i)^{[j]}$
for $i\in I$, forms an open cover
for $U^{[j]}$.
Given a subset $V\sub U$, let
$\lambda_V\colon V\to U$
and $\mu_V\colon V^{[j]}\to U^{[j]}$
be the inclusion maps.
Then, by {\bf\ref{compincov}}, the map $\tau_j\circ \Psi$
will be continuous if we can show
that $C^k(\mu_V,F)\circ \tau_j\circ \Psi$,
$f\mto \tau_j(\Psi(f))|_{V^{[j]}}$ is continuous
for each set~$V$ in our open cover of $U^{[j]}$.
For $V=U_i$, we have
$C^k(\mu_{U_i},F)\circ \tau_j\circ \Psi
=\tau_{j,i}\circ \psi_i
\circ C^k(\lambda_{U_i},F)_{SDS}$,
where $C^k(\lambda_{U_i},F)_{SDS}$
is continuous by Lemma~\ref{puback},
the mapping $\psi_i\colon C^k(U_i,F)_{SDS}\to C^k(U_i,F)_{BGN}$, $f\mto f$
is continuous by what has already been shown,
while the map $\tau_{j,i}\colon C^k(U_i,F)_{BGN}\to
C((U_i)^{[j]},F)$, $f\mto f^{[j]}$
is continuous by definition of the topology on
$C^k(U_i,F)_{BGN}$.
Thus $C^k(\mu_{U_i},F)\circ \tau_j\circ \Psi$
is continuous.
It only remains to see that also
$C^k(\mu_V,F)\circ \tau_j\circ \Psi$
is continuous with $V:=(U^{[j-1]})^{]1[}$.
For $f\in C^k(U,F)$ and $(x,y,t)\in
(U^{[j-1]})^{]1[}\sub U^{[j-1]}\times (\K^d)^{[j-1]}\times \K$,
we have
%
\begin{eqnarray}
\hspace*{-28mm}
\lefteqn{(C^k(\mu_V,F)\circ \tau_j\circ \Psi)(f)(x,y,t)}\qquad\quad\notag\\
&= & f^{[j]}(x,y,t)\; =\;
{\textstyle \frac{1}{t}}(f^{[j-1]}(x+ty)-f^{[j-1]}(x))\,.\,\,\label{prefinal}
\end{eqnarray}
The map $g\colon (U^{[j-1]})^{]1[}\to \K$,
$(x,y,t)\mto \frac{1}{t}$ is continuous,
and both of the maps
$s\colon (U^{[j-1]})^{]1[}\to U^{[j-1]}$,
$s(x,y,t):=x+ty$
and $p \colon (U^{[j-1]})^{]1[}\to U^{[j-1]}$,
$(x,y,t)\mto x$
are~$C^\infty$.
With the help of the multiplication operator\linebreak
$m_g\colon C(V,F)\to C(V,F)$ by~$g$,
we can rewrite~(\ref{prefinal}) in the form
\[
C^k(\mu_V,F)\circ \tau_j\circ \Psi
\;=\;
m_g\circ \big((C(s,F)\circ \tau_{j-1}\circ \Psi_{j-1})\,-\,
(C(p,F) \circ \tau_{j-1}\circ \Psi_{j-1})\big)\,.
\]
Here $\tau_{j-1}$ is continuous,
$\Psi_{j-1}$ is continuous by induction,
and all of the maps $m_g$,
$C(s,F)$ and $C(p,F)$ are continuous
by {\bf\ref{contfctsp}}.
Therefore the map $C^k(\mu_V,F)\circ \tau_j\circ \Psi$
is continuous, which completes the proof.\,\Punkt
\section{H\"{o}lder differentiable maps into
real locally convex spaces}\label{appreal}
%
%
For mappings from open sets
into real or complex locally convex spaces,
we provide a simpler characterization
of $C^{k,\sigma}$-maps now in terms
of the existence and H\"{o}lder continuity
of higher differentials.
They turn out to be the $C^{k,\sigma}$-maps
arising in the approach to $C^k$-maps
by Michal and Bastiani.
%
\begin{numba}\label{convapp}
Throughout this section, $\K\in \{\R,\C\}$,
$F$ is a locally convex topological $\K$-vector space,
$E$ a topological $\K$-vector space,
$U\sub E$ open,
and~$\sigma\in \;]0,1]$.
\end{numba}
%
%
\begin{defn}\label{defMB}
We say
that a map $f\colon U\to F$ is~$C^{0,\sigma}_{MB}$
if it is~$C^{0,\sigma}$.
The map $f\colon U\to F$ is $C^{1,\sigma}_{MB}$
if $d^0f:=f$ is~$C^{0,\sigma}$,
the directional (real, resp., complex) derivative
\[
d^1f(x,y)\, :=\, df(x,y)\, :=\, (D_yf)(x)\,:=\,
{\textstyle\frac{d}{dt}\big|_{t=0}}f(x+ty)
\]
exists for each $x\in U$ and $y\in E$,
and the map $df\colon U\times E\to F$
so obtained is~$C^{0,\sigma}$.
The map~$f$ is~$C^{2,\sigma}_{MB}$
if both~$f$ and~$df$ are~$C^{1,\sigma}_{MB}$;
we then define $d^2f:=d(df)\colon U\times E^3\to F$.
Recursively,
having defined $C^{k-1,\sigma}_{MB}$-maps
for some $k\in \N\cup\{\infty\}$,
we say that~$f$ is~$C^{k,\sigma}_{MB}$
if~$f$ is~$C^{1,\sigma}_{MB}$
and~$df$ is~$C^{k-1,\sigma}_{MB}$;
we then define $d^kf:=d^{k-1}(df)\colon U\times E^{2^k-1}\to F$.
\end{defn}
Replacing $C^{0,\sigma}$-maps
by continuous maps in the preceding definition,
we obtain the usual definition of
$C^k_{MB}$-maps
(also known as Keller's $C^k_c$-maps);
see~\cite{RES} and~\cite{GaN}
for further information.
Apparently, every $C^{k,\sigma}_{MB}$-map
is~$C^k_{MB}$.\\[2.5mm]
Our goal is the following result:
%
%
\begin{thm}\label{viadiff}
In the situation of~{\bf\ref{convapp}},
a map $f\colon U\to F$ is~$C^{k,\sigma}_{BGN}$
if and only if~$f$ is~$C^{k,\sigma}_{MB}$.
\end{thm}
Before we prove Theorem~\ref{viadiff},
let us mention an alternative description
of $C^{k,\sigma}_{MB}$-maps which
some readers may find more natural.
\begin{rem}
If $f$ is~$C^{k,\sigma}_{MB}$,
then the
iterated directional derivatives
\[
d^{(j)}f(x,v_1,\ldots, v_j)\; :=\;
(D_{v_1}\cdots D_{v_j}f)(x)
\]
exist for all $j\in \N$ such that $j\leq k$,
$x\in U$ and $v_1,\ldots, v_j\in E$.
Also, $d^{(j)}f$ is a partial map of~$d^jf$
(see \cite[Claim~2 on p.\,50]{RES}
or \cite[\S\,1]{GaN})
and hence~$C^{0,\sigma}$.
If, conversely, $d^{(j)}f$ exists for each
$j\in \N_0$ such that $j\leq k$
and is~$C^{0,\sigma}$,
then~$f$ is~$C^k$ (see \cite[Lemma~1.14]{RES}
or \cite[\S\,1]{GaN})
and formulas are available
which express~$d^jf$ in terms of
the maps~$d^{(i)}f$ with $i\leq j$
(see \cite[Claim~1 on p.\,49]{RES}
or \cite[\S\,1]{GaN}).
Inspecting these formulas,
we readily see that~$d^jf$ is~$C^{0,\sigma}$
for each $j\in\N_0$
such that $j\leq k$ and thus~$f$ is~$C^{k,\sigma}_{MB}$.
Although this second description
involving the higher differentials~$d^{(j)}f$
may look more palatable,
it is much more convenient for
our inductive proofs below
to use Definition~\ref{defMB}
and the (somewhat unfamiliar)
iterated differentials~$d^jf$.
\end{rem}
The proof of Theorem~\ref{viadiff}
is based on two lemmas concerning parameter-dependent
integrals. Recall that, given a continuous map
$\gamma\colon [0,1]\to F$ to a real or complex
locally convex space~$F$,
an element $z\in F$ is called
the \emph{weak integral} of~$\gamma$
if $\lambda(z)=\int_0^1 \lambda(\gamma(t))\, dt$
for all continuous linear
functionals $\lambda\in F'$.
Then~$z$ is uniquely determined,
and we write $\int_0^1\gamma(t)\,dt:=z$.
If~$F$ is complete (or at least sequentially complete),
then the weak integral $\int_0^1\gamma(t)\,dt$
always exists (see~\cite{GaN}).
%
%
%
\begin{la}\label{integhoel}
Let $k\in \N_0$, $W\sub \K$ be an open
subset such that $[0,1]\sub W$
and $h\colon U\times W \to F$
be a $C^{k,\sigma}_{MB}$-map.
Suppose that the weak integral
\[
g(x)\, :=\, \int_0^1 h(x,t)\; dt
\]
exists in~$F$ for each $x\in U$, and suppose that the
weak integrals
\[
\int_0^1 d_1^jh(x,y,t)\; dt
\]
exist for all $j\in \{1,\ldots, k\}$, $x\in U$
and $y\in E^{2^j-1}$,
where $d_1^jh(x,y,t)
:=d^j(h(\sbull, t))(y)$
for $x\in U$, $y\in E^{2^j-1}$ and $t\in W$
denotes the $j$-th iterated partial
differential of~$h$ with respect to the first variable.
Then $g\colon U\to F$ is a $C^{k,\sigma}_{MB}$-map, and
%
%
\begin{equation}\label{anotheq}
d^jg(x,y)\; =\;
\int_0^1 d_1^jh(x,y,t)\; dt
\end{equation}
for all $j\in \{1,\ldots, k\}$,
$x\in U$, and $y\in E^{2^j-1}$. 
\end{la}
\begin{proof}
The proof is by induction on~$k$.\\[2.5mm]
\emph{The case $k=0$}:
To see that~$g$ is~$C^{0,\sigma}$,
let~$q$ be a continuous seminorm on~$F$
and $x_0\in U$.
Given $t\in [0,1]$ and $r>0$,
set $B_r(t):=\{s\in [0,1]\colon
|s-t|<r\}$.
Since $h|_{U\times [0,1]}$ is~$C^{0,\sigma}$,
for each $t\in [0,1]$ we find
a neighbourhood~$V_t\sub U$ of~$x_0$,
$r_t>0$,
a gauge~$p_t$ on~$E$ and $C_t\in [0,\infty[$
such that
%
\begin{equation}\label{bizzz}
\|h(x,s)-h(y,r)\|_q
\; \leq \; \max\{ (\|x-y\|_{p_t})^\sigma,
C_t|s-r|^\sigma\}
\end{equation}
for all $x,y\in V_t$ and $r,s\in B_{2 r_t}(t)$
(compare Lemma~\ref{findfund}
and the proof of Lemma~\ref{prodhoel}).
Since~$h$ and~$q$ are continuous and $[0,1]$
is compact, after shrinking~$V_t$
we may assume that
\[
s_t\; :=\; \sup\{\|h(x,s)-h(y,r)\|_q\colon \mbox{$x,y\in V_t$,
$s,r\in [0,1]$}\}\; <\; \infty \, .
\]
There is a finite subset $\Phi\sub [0,1]$
such that $[0,1]\sub \bigcup_{t\in \Phi}B_{r_t}(t)$.
We define
$V:=\bigcap_{t\in \Phi}V_t$,
$p(x):=\max\{p_t(x)\colon t\in\Phi\}$
for $x\in E$
and let~$C$ be the maximum of
the numbers~$C_t$ and $\frac{s_t}{(r_t)^\sigma}$
for $t\in \Phi$.
Then
%
\begin{equation}\label{bizz2}
\|h(x,s)-h(y,r)\|_q
\; \leq \; \max\{ (\|x-y\|_p)^\sigma,
C|s-r|^\sigma\}
\end{equation}
for all $x,y\in V$ and $r,s\in [0,1]$.
In fact, there is $t\in [0,1]$
such that $r\in B_{r_t}(t)$.
If $s\in B_{2r_t}(t)$, then~(\ref{bizz2})
holds as a consequence of~(\ref{bizzz}).
Otherwise, $|s-r|\geq r_t$ and hence
$\|h(x,s)-h(y,r)\|_q\leq s_t=\frac{s_t}{(r_t)^\sigma} (r_t)^\sigma
\leq C (r_t)^\sigma\leq C |s-r|^\sigma$,
entailing that~(\ref{bizz2})
also holds in this case.
For all $x,y\in V$, using~(\ref{bizz2}) we obtain
\begin{eqnarray*}
\|g(x)-g(y)\|_q & = & \left\|\int_0^1(h(x,t)-h(y,t))\;dt\right\|_q
\; \leq \; \int_0^1\|h(x,t)-h(y,t)\|_q\;dt\\
& \leq & \int_0^1(\|x-y\|_p)^\sigma \;dt
\; =\; (\|x-y\|_p)^\sigma\, .
\end{eqnarray*}
Hence $g$ is $C^{0,\sigma}$, as required.\\[2.5mm]
\emph{Induction step.} Assume that $k\in \N$
and assume that the assertion of the lemma holds
if~$k$ is replaced with $k-1$.
By hypothesis,
$\phi(x,y):=
\int_0^1 d_1h(x,y,t)\,dt$ exists in~$F$
for all $x\in U$, $y\in E$,
and~$g$ as well as the map $\phi\colon U\times E\to F$
so obtained is~$C^{0,\sigma}$, by the case
$k=0$.
The map $d_1h \colon (U\times E)\times W \to F$
is $C^{k-1,\sigma}_{MB}$
and satisfies analogous conditions
as~$h$, with $k-1$ in place of~$k$,
because $d_1^j(d_1h)=d^{j+1}_1h$
for each $j\in \{0,1,\ldots, k-1\}$.
Hence~$\phi$ is $C^{k-1,\sigma}_{MB}$,
by induction.
We claim that the directional
derivative $dg(x,y)$
exists for all $x\in U$ and $y\in E$,
and is given by
$dg(x,y)=\phi(x,y)$.
Since~$\phi$ is $C^{k-1,\sigma}_{MB}$,
this entails that~$g$ is $C^{1,\sigma}_{MB}$
with $dg=\phi$ a $C^{k-1,\sigma}_{MB}$-map,
and thus~$g$ is
$C^{k,\sigma}_{MB}$.
Furthermore, $d^jg(x,y)=d^{j-1}(dg)(x,y)=d^{j-1}\phi(x,y)=
\int_0^1 d_1^{j-1}(d_1h)(x,y,t)\,dt
=\int_0^1 d_1^jh(x,y,t)\,dt$
for $x\in U$ and $y\in E^{2^j-1}$,
using the induction hypothesis
for the third equality.
Therefore~(\ref{anotheq}) holds.
To complete the inductive proof,
it remains to verify the claim.
To this end, fix $(x,y)\in U\times E$ and
pick $\ve >0$ such that
$x+B^\K_{2\ve }(0)y\sub U$.
Then
\[
R \colon B^\K_\ve(0)\times W
\times B^\K_2(0) \to F\, , \quad
R(s,t,r)\, :=\, d_1h(x+rsy, y,t)
\]
is a $C^{0,\sigma}$-map.
The weak integral
$K(s,t) :=\int_0^1 R(s,t,r)\,dr$
exists in~$F$ for all $(s,t)\in  B^\K_\ve(0)\times W$.
In fact, $K(s,t)=\frac{h(x+sy,t)-h(x,t)}{s}$
satisfies the defining property
of the weak integral if $s\not=0$,
by the Fundamental Theorem of Calculus~\cite[Theorem~1.5]{RES}.
For $s=0$, the integrand is constant
and hence
coincides with the weak integral;
thus $K(0,t)=d_1h(x,y,t)$.
By the case $k=0$, the map
$K\colon B^\K_\ve(0) \times W \to F$
is $C^{0,\sigma}$.
The weak integral
$L(s):=\int_0^1K(s,t)\,dt$\linebreak
exists
in~$F$ for each $s\in B^\K_\ve(0)$
since $L(s)=\frac{g(x+sy)-g(x)}{s}$
satisfies the defining property
of the weak integral if $s\not=0$,
while $L(0)=\int_0^1 d_1h(x,y,t)\, dt$
exists by hypothesis.
Hence $L\colon B^\K_\ve(0)\to F$
is $C^{0,\sigma}$ be the case $k=0$
and therefore continuous.
As a consequence,
$dg(x,y)=\lim_{s\to 0}\frac{g(x+sy)-g(x)}{s}
=\lim_{s\to 0}L(s)=L(0)=\int_0^1 d_1h(x,y,t)\, dt
=\phi(x,y)$ exists and is of the desired form.
\end{proof}
%
%
\begin{la}\label{convint}
Let $h\colon U\times W\to F$ be a $C^{1,\sigma}_{MB}$-map,
where
$W\sub \K$ is an open
neighbourhood of $[0,1]$.
Suppose that the weak integral
\[
g(x)\; :=\; \int_0^1 h(x,t)\; dt
\]
exists in~$F$ for each $x\in U$, and defines
a map $g\colon U\to F$ which is~$C^{1,\sigma}_{MB}$.
Then the weak integral $\int_0^1 d_1h(x,y,t)\, dt$
exists in~$F$ for all $x\in U$ and $y\in E$,
and it is given by
\[
\int_0^1 d_1h(x,y,t)\; dt \; =\;
dg(x,y)\,.
\]
\end{la}
\begin{proof}
We consider $g$ and $h$ as maps into the
completion~$\wt{F}$ of~$F$.
The weak integral in question exists in~$\wt{F}$.
By the preceding lemma, it coincides with $dg(x,y)$
and therefore is an element of~$F$.
\end{proof}
{\bf Proof of Theorem~\ref{viadiff}.}
We may assume that $k\in \N_0$. It is clear that
each $C^{k,\sigma}_{BGN}$-map is $C^{k,\sigma}_{MB}$
(noting that $d^jf$ is a partial map of~$f^{[j]}$, for each $j\leq k$).
The converse direction is proved by induction on~$k$.
If~$k=0$, then $C^{0,\sigma}_{BGN}$-maps
and $C^{0,\sigma}_{MB}$-maps coincide by definition.
Now suppose that~$f$ is a $C^{k,\sigma}_{MB}$-map,
where $k\in \N$, and suppose that
each $C^{k-1,\sigma}_{MB}$-map is~$C^{k-1,\sigma}_{BGN}$.
Then~$f$ is $C^{k-1,\sigma}_{BGN}$, by induction.
If we can show that the map
\[
g\colon U^{[1]}\to F\, , \quad
(x,y,t)\mto\left\{
\begin{array}{cl}
\frac{f(x+ty)-f(x)}{t}
& \;\mbox{if $\,t\in \K^\times$};\\
df(x,y) & \;\mbox{if $\, t=0$}
\end{array}
\right.
\]
is $C^{k-1,\sigma}_{BGN}$,
then~$f$ will be~$C^{1,\sigma}_{BGN}$
with~$f^{[1]}=g$ a $C^{k-1,\sigma}_{BGN}$-map,
and thus~$f$ will be $C^{k,\sigma}_{BGN}$, as required.
To see that~$g$ is $C^{k-1,\sigma}_{BGN}$,
note first that
$g|_{U^{]1[}}$
takes $(x,y,t)$ to $\frac{f(x+ty)-f(x)}{t}$.
Since~$f$ is $C^{k-1,\sigma}_{BGN}$,
it follows that $g|_{U^{]1[}}$ is $C^{k-1,\sigma}_{BGN}$.
Furthermore, by the Fundamental
Theorem of Calculus~\cite[Theorem~1.5]{RES},
we have
\[
\frac{f(x+ty)-f(x)}{t}\;=\;
\int_0^1 df(x+sty,y)\; ds
\]
for $(x,y,t)$ in the open neighbourhood
$V\!:=\! \{(x,y,t)\in U^{[1]}\! \colon
x\!+\!Kty \!\sub \!U\}$\linebreak
of $U\times E\times \{0\}$ in~$U^{[1]}$,
with $K\sub \K$ a compact neighbourhood
of $[0,1]$.~Thus
%
\begin{equation}\label{intform}
g(x,y,t)\;=\; \int_0^1df(x+sty,y)\; ds\; =\;
\int_0^1 h((x,y,t), \, s)\; ds
\end{equation}
for all $(x,y,t)\in V$,
where $h\colon V\times K^0 \to F$,
$h((x,y,t),\, s):=df(x+sty, y)$
is $C^{k-1,\sigma}_{BGN}$
since so is~$df$.
Set $d_1^{(j)}h(v;y_1,\ldots, y_j; s):=
(D_{y_1}\cdots D_{y_j}h(\sbull,s))(v)$
for $j\in \{1,\ldots, k-1\}$,
$v \in V$,
$y_1,\ldots, y_j\in E^{[1]}$
and $s\in K^0$.
In the proof of\linebreak
\cite[Proposition~7.4]{BGN},
it was shown that the weak integral
%
%
\begin{equation}\label{eqtem}
\int_0^1 d_1^{(j)}h(v;y_1,\ldots, y_j; s)\;ds
\end{equation}
exists in~$F$.
We now use that $d^j_1h(v,y,s)$
with $v\in V$, $y=(y_1,\ldots, y_{2^j-1})\in (E^{[1]})^{2^j-1}$,
$s\in K^0$
can be expressed in the form
\[
d^j_1h(v,y,s)\;=\;
\sum_{i=1}^j\sum_\ell n^{(i)}_\ell
d_1^{(i)}h(v;y_{\ell_1},\ldots, y_{\ell_i}; s)\,,
\]
where the second summation is over
all $\ell=(\ell_1,\ldots,\ell_i)\in \{1,\ldots, 2^j-1\}^i$,
and $n^{(i)}_\ell\in \N_0$ (cf.\
\cite[p.\,49, Claim~1]{RES}).
As a consequence, also the weak integrals
$\int_0^1 d_1^jh(v,y,s)\, ds$
exists for all $j\in \{1,\ldots, k-1\}$,
$v \in V$ and $y\in (E^{[1]})^{2^j-1}$
(since they are linear combinations
of those\vspace{1mm} in~(\ref{eqtem})).
Thus $h$ satisfies the hypotheses of Lemma~\ref{integhoel}
with $k-1$ in place of~$k$,
and hence~$g$ is $C^{k-1,\sigma}_{MB}$
and thus $C^{k-1,\sigma}_{BGN}$, by induction.
This completes the proof.\,\,\Punkt
{\footnotesize

{\bf Helge Gl\"{o}ckner}, TU Darmstadt, Fachbereich Mathematik AG~5,
Schlossgartenstr.\,7,\\
64289 Darmstadt, Germany. E-Mail: gloeckner\at{}mathematik.tu-darmstadt.de}
\end{document}